\documentclass[12pt]{elsarticle}
\usepackage{amsfonts,enumerate,amsmath,amssymb,mathrsfs,amsbsy}
\usepackage[colorlinks=true]{hyperref}

\usepackage{caption}
\biboptions{numbers,sort&compress}

\usepackage{xcolor}

\def\qed{\strut\hfill $\Box$}
\newtheorem{thm}{Theorem}[section]

\newtheorem{lem}[thm]{Lemma}

\newtheorem{rem}[thm]{Remark}

\newcommand{\thmref}[1]{Theorem~{\rm \ref{#1}}}
\newcommand{\lemref}[1]{Lemma~{\rm \ref{#1}}}

\def\para#1{\vskip .4\baselineskip\noindent{\bf #1}}
\numberwithin{equation}{section}
\allowdisplaybreaks \allowdisplaybreaks[4]
\begin{document}
	\begin{frontmatter}	
		\title{Stochastic averaging for the non-autonomous mixed stochastic differential equations with locally Lipschitz coefficients}
		
		\author[mymainaddress]{Ruifang Wang}
		\ead{wrfjy@yahoo.com}
		
		\author[mymainaddress,myfivearyaddress]{Yong Xu\corref{mycorrespondingauthor}}
		\cortext[mycorrespondingauthor]{Corresponding author}
		\ead{hsux3@nwpu.edu.cn}
		
		\author[mymainaddress]{Hongge Yue}
		\ead{yuehongge803@163.com}

		%	\author[mysecondaryaddress]{Bin Pei}
		%	\ead{binpei@hotmail.com}
		
		%	\author[myfivearyaddress]{Yuzhen Bai}
		%	%		\ead{xgdsxxpeibin@126.com}
		
		%	\author[mythirdaryaddress,myfouraryaddress]{J\"{u}rgen Kurths}
		%	\ead{Juergen.Kurths@pik-potsdam.de}
		
		\address[mymainaddress]{School of Mathematics and Statistics, Northwestern Polytechnical University, Xi'an, 710072, China}
		\address[myfivearyaddress]{MIIT Key Laboratory of Dynamics and Control of Complex Systems, Northwestern Polytechnical University, Xi'an, 710072, China	}
		%	\address[mysecondaryaddress]{School of Mathematical Sciences, Fudan University, Shanghai,  200433, China}
		%	
		%	\address[mythirdaryaddress]{Potsdam Institute for Climate Impact Research, Potsdam, 14412, Germany}
		%	\address[myfouraryaddress]{Department of Physics, Humboldt University Berlin, Berlin, 12489, Germany}

		\begin{abstract}
			This paper investigates a non-autonomous slow-fast system, which is generalized by 
			stochastic differential equations (SDEs) with locally  Lipschitz coefficients, subjected to standard Brownian motion (Bm) and fractional Brownian motion (fBm) with Hurst parameter $ 1/2<H<1. $ We concentrate on how to handle both types of integrals with respect to Bm and fBm and the locally Lispchitz continuity. The pathwise approach and the It\^{o} stochastic calculus are combined with the technique of stopping time  to establish the averaging principle where the averaged equation is defined. Then, the slow component of the original slow-fast system converges to the solution of the proposed averaged equation in the mean square sense is verified.
			%		70K70, 60H15, 34K33, 37B55, 60J75 
		\end{abstract}		
	\end{frontmatter}
	
	\section{Introduction}\label{sec-1}
	In this paper, we study the following SDEs driven by fBm and 
	standard Bm: 
	\begin{align}\label{orginal1}
	\left\{ \begin{array}{l}
	du_{t}^{\epsilon}=b_1\left( t,u_{t}^{\epsilon},v_{t}^{\epsilon} \right) dt+f_1\left( t,u_{t}^{\epsilon} \right) dW_{t}^{1}+g_1\left( t,u_{t}^{\epsilon} \right) dB_{t}^{H},\quad u_{0}^{\epsilon}=x\in \mathbb{R}^n,\\
	dv_{t}^{\epsilon}=\frac{1}{\epsilon}b_2\left( t,u_{t}^{\epsilon},v_{t}^{\epsilon} \right) dt+\frac{1}{\sqrt{\epsilon}}f_2\left( t,u_{t}^{\epsilon},v_{t}^{\epsilon} \right) dW_{t}^{2},\qquad\qquad\quad\  v_{0}^{\epsilon}=y\in \mathbb{R}^m,\\
	\end{array} \right.	
	\end{align}
	where $ t\in \left[ 0,T\right]  $ and $ \epsilon \in \left( 0,1 \right] $  is a small positive parameter which represents the ratio of the natural time scale between the slow variable $ u_{t}^{\epsilon}\in \mathbb{R}^n
	$  and fast variable $ v_{t}^{\epsilon}\in \mathbb{R}^m
	. $ Moreover, $ B^H=\left\{ B_{t}^{H},t\in \left[ 0,T \right] \right\} \left( H\in \left( 1/2,1 \right) \right) 
	$ is  $ d_1$-dimensional fBm,  $ W^1=\left\lbrace  W_{t}^{1}, \right. $   $\left. t\in \left[ 0,T \right] \right\rbrace  
	$ and $ W^2=\left\{ W_{t}^{2},\ t\in \left[ 0,T \right] \right\} 
	$  are $ d_2 $  and $ d_3$-dimensional standard Bm, respectively. Assume that $ W^1,\ W^2
	$ and $ B^H $  are mutually independent processes, and initial variables $ x,y $ are fixed and independent of $  \left( W^1,\ W^2,\ B^H \right) 
	. $ 
	
	Recall that fBm with Hurst index $ H\in \left(0,1 \right)  $ is a zero mean Gaussian process $ \left\lbrace B_t^H, t\geq 0\right\rbrace   $ with covariance function
	\begin{eqnarray}\label{en38}
	R_{H}(s, t)=\frac{1}{2}\left(s^{2 H}+t^{2 H}-|t-s|^{2 H}\right).
	\end{eqnarray}
	Notice that $ B^H $ is a standard Bm if $ H=1/2, $ but if $ H\neq 1/2, $ it does not have independent increments. Moreover, from (\ref{en38}) we deduce that, $ \mathbb{E}\left| B_t-B_s\right|^2=\left|  t- s\right|^{2H}.  $ As a
	consequence, the process $ B^H $  has $ \alpha $-H\"{o}lder continuous paths for all $ \alpha\in \left( 0,H\right).  $  
	Kolmogorov introduced the process $ B^H $  \cite{Kolmogorov1940Wienershe} in 1940, and Mandelbrot and Van Ness \cite{Mandelbrot1968Fractional} later named it as fBm in 1968.

	The averaging principle is a kind of approximate theorems to simplify stochastic dynamical systems with different time-scales, and the first related result about stochastic case was studied by Khasminskii \cite{khas1968on} in 1968. Since then, the averaging principle has been investigated by a number of investigators. For instance, in the case of the autonomous system: Givon \cite{givon2007strong}, Freidlin and Wentzell \cite{freidlin2012random},  Duan \cite{duan2014effective}, Xu and his co-workers \cite {xu2011averaging,xu2017stochastic} studied the averaging principle of SDEs. In addition,  
	the averaging principle of stochastic partial differential equations (SPDEs) also have been investigated by Cerrai \cite{cerrai2009khasminskii,cerrai2011averaging}, Wang and Roberts \cite{wang2012average}, Pei and Xu \cite{pei2017two,  pei2017averaging}  and others in recent years. For the non-autonomous case, Cerrai \cite{cerrai2017averaging} and Liu \cite{liu2020averaging} studied the averaging principle for non-autonomous slow-fast systems driven by Brownian motion, which is generalized by SPDEs and SDEs respectively. Moreover, Xu \cite{Xu2018Averaging} also studied a class of non-autonomous slow-fast systems of SPDEs driven by Gaussian noises and Poisson random measures, and proved the averaging principle in the sense of probability. However, either Gaussian noises or Poisson random measures can not describe the disturbances with long-range dependence. This paper is to consider a class of non-autonomous slow-fast system of SDEs driven by  standard Bm and fBm, and the convergence of averaging principle in the mean square sense is to be proved.

	First, the existence and uniqueness of  solutions for (\ref{orginal1}) is studied.
	Comparing with the work of Guerra and Nualart \cite{Guerra2008Stochastic},  
	the conditional expectation (given $ x, y $ and $ B^H $) in \cite{Guerra2008Stochastic} is extended to general expectation on 
	$\big( \varOmega ,\mathcal{F},\left\lbrace \mathcal{F}_t\right\rbrace _{t\geq0},\mathbb{P} \big) $ (where $ \left\lbrace \mathcal{F}_t\right\rbrace _{t\geq0} $ is the $ \sigma $-field generated by the random variables $  W_t^1, W_t^2, B^H_t $ and the $ \mathbb{P} $-null sets) in our work and the coefficients are assumed to satisfy locally Lipschitz conditions.   
	To overcome the problem of $ B^H $ is not a semimartingale if $ H \neq 1/2, $ we interpret the integral   $ \int_0^t{f_1\left( r,u_{r}^{\epsilon} \right)}dW_{r}^{1},\ \int_0^t{f_2\left( r,u_{r}^{\epsilon},v_{r}^{\epsilon} \right)}dW_{r}^{2}
	$ as an It\^{o} stochastic integral and the integral $ \int_0^t{g_1\left( r,u_{r}^{\epsilon} \right)}dB_{r}^{H}
	$  as a generalised Riemann-Stieltjes integral in the sense of Z\"{a}hle \cite{Zahle1998Integration}  in our work. 
	Then, the  It\^{o} stochastic calculus and the pathwise approach are combined to handle these two kinds of integrals and the Garsia-Rodemich-Rumsey inequality is used to overcome the problem of $ B^H $ will produce some higher order terms.
	Moreover, the technique of stopping time  is also used to deal with the problem of expectation and locally Lispchitz continuity.    
	
	Then, consider that the coefficients in equation (\ref{orginal1}) depend on time, the corresponding equation associated to the fast equation by fixed $ s>0  $ and frozen slow component $ x\in \mathbb{R}^n  $
	\begin{align}\label{en12}
	dv_t=b_2\left( s,x,v_t \right) dt+f_2\left( s,x,v_t \right) dW_{t}^{2},\quad v_0=y\in \mathbb{R}^m 
	\end{align} 
	is introduced. Further, the existence of the unique invariant measure  $ \mu ^{s,x}   $ for the equation (\ref{en12}) is proved and the averaged coefficient can be defined as 
	$$ \bar{b}_1\left( s,x \right) =\int_{\mathbb{R}^m}{b_1\left( s,x,z \right) \mu ^{s,x}\left( dz \right)}.  $$
	
	Finally, the detailed proof of the convergence result is presented by using the technique of time discretization and truncation. That is, under some reasonable assumptions, the purpose of this paper is to show the convergence of averaging principle in the mean square sense:
	\begin{align} 
	\underset{\epsilon \rightarrow 0}{\lim}\underset{t\in \left[ 0,T \right]}{\sup}\mathbb{E} \left\|  u_{t}^{\epsilon}-\bar{u}_t \right\|  _{\alpha}^{2}  =0, 
	\end{align}
	where $ \bar{u}_t $ is the solution to the so-called averaged equation as:
	\begin{align}\label{en45}
	d\bar{u}_t=\bar{b}_1\left( t,\bar{u}_t \right) dt+f_1\left( t,\bar{u}_t \right) dW_{t}^{1}+g_1\left( t,\bar{u}_t \right) dB_{t}^{H},\ \ \bar{u}_0=x\in \mathbb{R}^n.
	\end{align}

	This paper is organized as follows,
	%The structure of this work is as follows. 
	in Section \ref{sec-2}, some notations and assumptions which will be used in the analysis
	of equation (\ref{orginal1}) is introduced and the main results is claimed. Section \ref{sec-3} is committed to proving the existence and uniqueness of solutions to \ref{orginal1}. In Section \ref{sec-4}, the averaging principle is obtained by using the generalized Khasminskii method where the averaged equation is defined. Note that, $ C>0 $ with or without subscripts represents a general constant, the value of which may vary for different cases in this paper.

	\section{Preliminaries, assumptions and main result}\label{sec-2}
	Now, we recall some definitions and results that will be used throughout the paper.
	Let $ \left| \cdot \right|
	$  be the Euclidean norm, $ \left< \cdot ,\cdot \right> 
	$  be the Euclidean inner product and $ \left\|  \cdot \right\|  
	$  be the matrix norm. 
	
	Let $ 1/2<H<1, 1-H<\alpha<1/2 $ and $ k,l\in N^+, $ denote 
	$ W_0^{\alpha,\infty} (0,T;\mathbb{R}^k  )   $ the space of measurable functions $ f:[0,T]\rightarrow\mathbb{R}^k  $ such that 
	\begin{align}
	\|f\|_{\alpha, \infty}:=\sup _{t \in[0, T]}\left\| f\left( t\right) \right\|_\alpha <\infty,\nonumber
	\end{align}
	where
	\begin{align}
	\left\| f\left( t\right) \right\|_\alpha:= |f(t)|+\int_{0}^{t} \frac{|f(t)-f(s)|}{(t-s)^{\alpha+1}} d s.\nonumber
	\end{align}
	
	For $ 0<\eta\leq 1,  $ let $ C^\eta (0,T;\mathbb{R}^k  )  $ be space of $ \eta $-H\"{o}lder continuous functions $ f:[0,T]\rightarrow\mathbb{R}^k,  $  equipped with the norm
	\begin{align}
	\|f\|_{\eta}:=\sup _{t \in[0, T]}|f(t)|+\sup _{0 \leq s<t \leq T} \frac{|f(t)-f(s)|}{(t-s)^{\eta}}<\infty.\nonumber
	\end{align} 
	
	%  Fix the parameter $ \alpha  $ such that $  $  
	Denote by $ W_T^{1-\alpha,\infty}(0,T;\mathbb{R}^k  ) $ the space of measurable
	functions $ g : [0, T] \rightarrow \mathbb{R}^k $ such that
	\begin{align}
	\|g\|_{1-\alpha, \infty, T}:=\sup _{0<s<t<T}\Big(\frac{|g(t)-g(s)|}{(t-s)^{1-\alpha}}+\int_{s}^{t} \frac{|g(r)-g(s)|}{(r-s)^{2-\alpha}} d r\Big)<\infty.\nonumber
	\end{align}
	Moreover, denote by $ W_{0}^{\alpha, 1}(0,T;\mathbb{R}^l  )$ the space of measurable functions $ f:[0, T] \rightarrow \mathbb{R}^{l} $ such that
	\begin{align}
	\|f\|_{\alpha, 1}:=\int_{0}^{T} \frac{|f(s)|}{s^{\alpha}} ds+\int_{0}^{T} \int_{0}^{s} \frac{|f(s)-f(r)|}{(s-r)^{\alpha+1}} d r ds<\infty.\nonumber
	\end{align}
	Then, if $ f \in W_{0}^{\alpha, 1}(0,T;\mathbb{R}^l  ) $   and  $  g \in W_{T}^{1-\alpha, \infty}(0,T;\mathbb{R}^k  ),  $ for any $ t\in [0,T], $  we know that $ \int_{0}^{t} fdg  $ exists, and have   
	\begin{align}
	\Big|\int_{0}^{t} f \mathrm{d} g\Big| \leq \Lambda_{\alpha}^{0,t}(g)\|f\|_{\alpha, 1},\nonumber
	\end{align}
	where 
	\begin{align}
	\Lambda_{\alpha}^{0,T}(g) :=\frac{1}{\Gamma(1-\alpha)} \sup _{0<s<t<T}\left|\left(D_{t-}^{1-\alpha} g_{t-}\right)(s)\right|  
	\leq \frac{1}{\Gamma(1-\alpha) \Gamma(\alpha)}\|g\|_{1-\alpha, \infty, T}
	<\infty,\nonumber
	\end{align}
	and $ \left(D_{t-}^{1-\alpha} g_{t-}\right)(s) $ is the Weyl derivatives \cite{Samko1993Fractional} of $ g $.
	\begin{rem}\label{rem2.5}
		In particular, the trajectories of fBm $ B^H  \left(  H > 1/2 \right)  $ belong to $ W_T^{1-\alpha,\infty}(0,T;\mathbb{R}^k  )$ where $\left( 1-H  <\alpha<1/2\right) . $ As a consequence,  if the trajectories of a stochastic process  $ u=\left\lbrace u_t, t\in [0,T]\right\rbrace  $  belong to the space $ W_0^{\alpha,1}(0,T;\mathbb{R}^l  ), $ the generalised Riemann-Stieltjes
		integrals $ \int_{0}^{T} u_{s} d B_{s}^{H} $ exists, and the following estimate holds
		\begin{align}\label{en3}
		\Big|\int_{0}^{T} u_{s} d B_{s}^H\Big| \leq \Lambda_{\alpha}^{0,T}\left(B^{H}\right)\|u\|_{\alpha, 1}, 
		\end{align}
		where $ \Lambda_{\alpha}^{0,T}\left(B^{H}\right)
		:=\frac{1}{\Gamma(1-\alpha)} \sup _{0<s<t<T}\left|\left(D_{t-}^{1-\alpha} B^H_{t-}\right)(s)\right|  
		$ has moments of all orders (see Lemma 7.5 in Nualart and R{\u{a}}{\c{s}}canu \cite{Nualart2002Differential}).  
	\end{rem}

	The following lemma is the so-called Garsia-Rodemich-Rumsey inequality (see Theorem 1.4 in \cite{Garsia1974Monotonicity}):  
	\begin{lem}
		For any $ p\geq 1 $ and $ \theta>p^{-1}, $ there exists some constant $ C_{\theta,p} > 0 $ such that for any continuous function $ f $ on $\left[0,T \right],  $ have
		\begin{align}\label{en41}
		|f(t)-f(s)|^{p} \leq C_{\theta, p}|t-s|^{\theta p-1} \int_{0}^{T} \int_{0}^{T} \frac{|f(x)-f(y)|^{p}}{|x-y|^{\theta p+1}} dx dy.
		\end{align}
	\end{lem}  
	
	%Now, given a complete probability space $ \left( \varOmega ,\mathcal{F},\mathbb{P} \right) 
	%$  and fix a time interval $  [0,	T]. $ 
	%Consider the following bigger filtrations $  \left\{ \mathcal{G}_t \right\} _{t\in \left[ 0,T \right]}
	%: $
	%\begin{enumerate}[1.]
	%	\item 	$ \left\{ \mathcal{G}_t \right\} _{t\in \left[ 0,T \right]}
	%	$  is right-continuous and $ \mathcal{G}_0
	%	$  contains the  $ \mathbb{P} $-null sets.
	%	\item  $ x $  and $ B^H $  are  $ \mathcal{G}_0
	%	$-measurable, and $ W $  is a  $ \mathcal{G}_t
	%	$-Brownian motion.
	%\end{enumerate}
	%Denote $ \hat{\mathcal{F}}_t
	%$  is the  $ \sigma $-field generated by $ \left\{ x,\ B^H,\ W_s,\ s\in \left[ 0,t \right] \right\} 
	%$  and the  $ \mathbb{P} $-null sets, so we can get that $  \hat{\mathcal{F}}_t\subset \mathcal{G}_t
	%. $ Moreover, we denote by  $ \mathbb{E} $ the conditional expectation given $  \hat{\mathcal{F}}_0
	%, $ that is, given  $ x $ and $  B^H. $
	In this paper, the following maps 
	\begin{align}
	&b_1:\left[0, \infty \right) \times \mathbb{R}^n\times \mathbb{R}^m\times \Omega \rightarrow \mathbb{R}^n;\cr
	&f_1:\left[0, \infty \right)\times \mathbb{R}^n\times \Omega \rightarrow \mathbb{R}^{n\times d_2};\cr
	&g_1:\left[0, \infty \right)\times \mathbb{R}^n\times \Omega \rightarrow \mathbb{R}^{n\times d_1};\cr
	&b_2:\left[ 0,\infty\right)\times \mathbb{R}^n\times \mathbb{R}^m\times \varOmega \rightarrow \mathbb{R}^m;\cr &f_2:\left[0, \infty \right)\times \mathbb{R}^n\times \mathbb{R}^m\times \varOmega \rightarrow \mathbb{R}^{m\times d_3} \nonumber
	\end{align}
	are continuous. Then, the following assumptions are supposed to hold for  $ \mathbb{P} $-almost all $  \omega \in \varOmega: $ 
	
	\begin{enumerate}[({A}1)]
		\item 
		(a) For any $ R\in\mathbb{R}, $ $ y\in \mathbb{R}^m $ and $ x_i\in \mathbb{R}^n (i=1,2)   $ with  $ \left|x_i \right|\leq R, $ there exist some constants  $  \theta_1 \geq 0, $  such that
		$$
		\left| b_1\left( t,x_1,y  \right) -b_1\left( t,x_2,y  \right) \right|+\left\|  f_1\left( t,x_1 \right) -f_1\left( t,x_2 \right) \right\|  \leq C_{R,T}\big( 1+\left| y\right| ^{\theta_1} \big)\left| x_1-x_2 \right|.
		$$
		(b) For any  $ x\in \mathbb{R}^n $ and $ y_1,y_2\in \mathbb{R}^m,   $ there exist some constants  $  \theta_2,\theta_3 \geq 0  $ and $ 0<\kappa\leq 1, $ such that
		$$
		\left| b_1\left( t,x ,y_1 \right) -b_1\left( t,x,y_2 \right) \right| \leq C_{T}  |y_1-y_2| \big(1+\left| x\right|^{\theta_2} +\left| y_1\right|^{\theta_3}+\left| y_2\right|^{\theta_3}   \big);
		$$
		$$
		\left| b_1\left( t,x_1,y_1 \right) -b_1\left( s,x_1,y_1 \right) \right|+\left\|  f_1\left( t,x_1 \right) -f_1\left( s,x_1 \right) \right\|  \leq C_{T}  |t-s|^{\kappa}\big(1+\left| x_1\right|^{\theta_2} +\left| y\right|^{\theta_3}  \big);
		$$
		$$
		\left| b_1\left( t,x_1,y_1 \right) \right|+\left\|  f_1\left( t,x_1 \right) \right\|  \leq C_{T}\left( 1+\left| x_1 \right|+\left| y_1 \right| \right) .
		$$
		\item (a) The mapping $ g_1 $  is  continuously differentiable in  $ x\in \mathbb{R}^n. $ For any $ R\in\mathbb{R} $ and $ x_i\in \mathbb{R}^n (i=1,2)  $ with  $ \left|x_i \right|\leq R, $ there exist some constants  $  0< \gamma \le 1, $  such that
		$$
		\left\|  g_1\left( t,x_1 \right) -g_1\left( t,x_2 \right) \right\|  \leq C_{R,T}\left| x_1-x_2 \right|;
		$$
		$$
		\left\|  \nabla _{x_1}g_1\left( t,x_1 \right) -\nabla _{x_2}g_1\left( t,x_2 \right) \right\|   \leq C_{R,T}\left| x_1-x_2 \right|^{\gamma},
		$$
		where $ \triangledown _x
		$  is the standard gradient with respect to the variable $  x. $\\
		(b)  For any $ x \in \mathbb{R}^n $ and $
		s,t\in \left[ 0,T \right],   $ there exist some constants  $  0<\beta  \le 1, $  such that
		$$
		\left\|  g_1\left( t,x  \right) \right\|  \le C_{T} \left( 1+ \left|  x  \right|  \right);
		$$
		$$
		\left\|  \nabla _{x}g_1\left( t,x \right) -\nabla _{x}g_1\left( s,x \right) \right\|  +\left\|  g_1\left( t,x \right) -g_1\left( s,x \right) \right\|  \leq C_{T}|t-s|^{\beta}.
		$$
		\item 
		(a) 
		For any $ x \in \mathbb{R}^n  $ and $ t\in \left[ 0,T \right],  $ the mapping $ b_2\left( t,x,\cdot\right)
		%	:\mathbb{R}^n\times \mathbb{R}^m\rightarrow \mathbb{R}^m   
		$  is  locally  Lipschitz continuous and  $ f_2\left( t,\cdot,\cdot\right)
		%:	\mathbb{R}^n\times \mathbb{R}^m\rightarrow \mathbb{R}^{m\times d_3} 
		$  is  Lipschitz continuous. \\ (b) For any $ x_i\in \mathbb{R}^n,\ y_i \in \mathbb{R}^m $ and $ t\in \left[ 0,T \right],  $ there exist some constants $  \alpha_1,\alpha_2 >0  $ and $ \iota\in\left( 0,1\right] , $ such that
		$$
		\left| b_2\left( t,x_1,y_1 \right) -b_2\left( t,x_2,y_1 \right) \right| \leq C_{T}   \left| x_1-x_2 \right| \big(1+\left| x_1\right|^{\alpha_1} +\left| x_2\right|^{\alpha_1}+\left| y_1\right|^{\alpha_2}  \big); 
		$$ 
		$$
		\left| b_2\left( t,x_1,y_1 \right) -b_2\left( s,x_1,y_1 \right) \right|+\left| f_2\left( t,x_1,y_1 \right) -f_2\left( s,x_1,y_1 \right) \right| \leq C_{T}  \left| t-s \right|^{\iota}\big(1+\left| x_1\right|^{\alpha_1} +\left| y_1\right|^{\alpha_2}  \big); 
		$$
		$$
		\left| b_2\left( t,x_1,y_1 \right) \right|+\left\|  f_2\left( t,x_1,y_1 \right) \right\| \leq C_{T}\left( 1+\left| x_1 \right|+\left| y_1 \right| \right) .
		$$
		\item Assume that for any $t\in \mathbb{R},\ x\in \mathbb{R}^n $ and $ y\in \mathbb{R}^m, $  $  b_1(t,x,y) $ and $ f_1(t,x) $ are bounded .
		%	 $  f_1 $ is bounded: exists some $ K > 0, $ such that $ \sup _{t \in[0, T],x\in \mathbb{R}^n} \left| f_1\left(t,x \right)\right|  \leq K. $
		\item (Strict monotonicity condition:) For any $ t\in \left[ 0,+\infty \right) ,\ x\in \mathbb{R}^n,\ y_1,y_2\in \mathbb{R}^m, $  there exist constants $ \beta _1 >0, $  such that
		\begin{align}
		2\left< y_1-y_2,b_2\left( t,x,y_1 \right) -b_2\left( t,x,y_2 \right) \right> +\left\|  f_2\left( t,x,y_1 \right) -f_2\left( t,x,y_2 \right) \right\|  ^2\le -\beta _1\left| y_1-y_2 \right|^2.\nonumber 
		\end{align}
		(Strict coercivity condition:) For some fixed $ p\geq 2, $  any $ t\in \left[ 0,+\infty \right) ,\ x\in \mathbb{R}^n  $ and $y\in \mathbb{R}^m, $  there exist constants $ C_{p,T},\ \beta_p>0, $ such that
		\begin{eqnarray}
		2\left< y ,b_2\left( t,x,y  \right) \right> +(p-1)\left\|  f_2\left( t,x,y  \right) \right\|  ^2\le -\beta _p\left| y  \right|^2+C_{p,T} ( 1+\left| x \right|^2 ).\nonumber 
		\end{eqnarray}
	\end{enumerate}
	\begin{rem}
		The existence and uniqueness of solutions for original equation (\ref{orginal1}) are guaranteed by the assumptions (A1)-(A4). Strict monotonicity condition guarantee the exponential ergodicity (see Lemma \ref{lem5.2} in Section \ref{sec-4}) holds and strict coercivity condition is used to ensures the existence of invariant measures for the frozen equation  (see \lemref{lem4.1} in Section \ref{sec-4}). 
		%关于A5的一些解释	
	\end{rem}

	Under the above assumptions, the main
	result of this paper is claimed as follows: 
	\begin{thm}\label{thm3.1}
		Assume that the conditions (A1)-(A4) hold. Then, for any $\alpha \in \left(  1-H, \frac{1}{2} \land\beta\right.  $  $\left.  \land  \frac{\gamma}{2} \right)    $ and $  \epsilon \in \left( 0,1\right],  $  
		there exists a unique solution  $ ( u_t^{\epsilon },v_t^{\epsilon } )\in (W_0^{\alpha,\infty} (0,T;\mathbb{R}^n  ), W_0^{\alpha,\infty} (0,T;\mathbb{R}^m  ) ) $ of equation (\ref{orginal1}).
%		\begin{align}
%		\left\{ \begin{array}{l}
%		u_{t}^{\epsilon }=x+\int _0^t{b_1\left( r,u_{r}^{\epsilon },v_{r}^{\epsilon } \right) dr}+\int_0 ^t{f_1\left( r,u_{r}^{\epsilon } \right) dW_{r}^{1}}+\int_0 ^t{g_1\left( r,u_{r}^{\epsilon } \right) dB_{r}^{H}},\\
%		v_{t}^{\epsilon }=y+\frac{1}{\epsilon }\int_0 ^t{b_2\left( r,u_{r}^{\epsilon },v_{r}^{\epsilon } \right) dr}+\frac{1}{\sqrt{\epsilon }}\int_0 ^t{f_2\left( r,v_{r}^{\epsilon } \right) dW_{r}^{2}}.\\
%		\end{array} \right. 
%		\end{align}
	\end{thm}
	\begin{thm}\label{thm2.1}
		Assume that the conditions (A1)-(A5) hold. Then, for any $ \epsilon \in \left( 0,1\right]    $ and $ \alpha \in \left(  1-H, \frac{1}{2} \land\beta  \land  \frac{\gamma}{2} \right),  $
		we have
		\begin{align}\label{en35}
		\underset{\epsilon \rightarrow 0}{\lim}\underset{t\in \left[ 0,T \right]}{\sup}\mathbb{E} \left\|  u_{t}^{\epsilon}-\bar{u}_t \right\|  _{\alpha}^{2}  =0, 
		\end{align}
		where $ \bar{u}_t $ is the solution of the corresponding averaged equation (\ref{en45}). 
	\end{thm}
	
	\section{Existence, uniqueness of the solutions}\label{sec-3}
	In this section, we aim to prove the existence and uniqueness of solutions for a class of mixed SDEs (\ref{orginal1}) driven by fBm and Bm (\thmref{thm3.1}). An auxiliary equation is introduced at first. Then, some estimates for this auxiliary equation are presented. Finally, the existence and uniqueness of solutions for original equation (\ref{orginal1}) is proved by defining the stopping time.
	\subsection{Some a-priori estimates of $  (u_t^{\epsilon,k},v_t^{\epsilon,k} ) $}
	For any $ k\in \mathbb{N}, $ we define the following stopping time 
	\begin{align}
	\tau _k:=\text{inf}\left\{ t\geq 0:\varLambda _{\alpha}^{0,t}\left( B^H \right) \geq k \right\},
	\end{align}
	and study the following equation:
	\begin{align}\label{orginal4}
	\left\{ \begin{array}{l}
	du_{t}^{\epsilon}=b_{1,k}\left( t,u_{t}^{\epsilon},v_{t}^{\epsilon} \right) dt+f_{1,k}\left( t,u_{t}^{\epsilon} \right) dW_{t}^{1}+g_{1,k} \left( t,u_{t}^{\epsilon} \right) dB_{t}^{H,k},\quad u_{0}^{\epsilon}=x\in \mathbb{R}^n,\\
	dv_{t}^{\epsilon}=\frac{1}{\epsilon}b_{2,k}\left( t,u_{t}^{\epsilon},v_{t}^{\epsilon} \right) dt+\frac{1}{\sqrt{\epsilon}}f_2\left( t,u_{t}^{\epsilon},v_{t}^{\epsilon} \right) dW_{t}^{2},\qquad\qquad\quad\quad \quad  v_{0}^{\epsilon}=y\in \mathbb{R}^m,\\
	\end{array} \right.	
	\end{align}
	where $ B^{H,k}_t =B^{H}_{t\land \tau _k}  $ and
	\begin{align}
	b_{i,k}\left( t,x,y \right) =\left\{ \begin{array}{c}
	b_i\left( t,x,y \right) ,\qquad \ \ \ \ \ \ \ \ \ \left| x \right|\leq k \ \text{and} \ \left| y \right|\leq k,\\
	b_i\left( t,x,yk/\left| y \right| \right) ,\qquad  \ \ \left| x \right|\leq k \ \text{and}\ \left| y \right|>k,\\
	b_i\left( t,xk/\left| x \right|,y \right) ,\qquad  \ \ \left| x \right|>k \ \text{and}\ \left| y \right|\leq k,\\
	b_i\left( t,xk/\left| x \right|,yk/\left| y \right| \right),  \   \left| x \right|>k \ \text{and}\ \left| y \right|>k,\\
	\end{array} \right. \nonumber
	\end{align}
	and 
	\begin{align}
	f_{1,k}\left( t,x \right) =\left\{ \begin{array}{c}
	f_1\left( t,x \right) ,\qquad\ \ \left| x \right|\leq k,\\
	f_1\left( t,xk/\left| x \right| \right) ,\  \left| x \right|>k,\\
	\end{array} \right. \qquad g_{1,k}\left( t,x \right) =\left\{ \begin{array}{c}
	g_1\left( t,x \right) ,\qquad\ \ \left| x \right|\leq k,\\
	g_1\left( t,xk/\left| x \right| \right) ,\  \left| x \right|>k.\\
	\end{array} \right.\nonumber
	\end{align}
	It is easy to know that the mapping $ b_{i,k}\left(  t,\cdot,\cdot \right),\ f_{1,k}\left(  t,\cdot  \right) $ and $ g_{1,k}\left(  t,\cdot \right) $ are Lipschitz continuous and satisfy all conditions in (A1)-(A3). Moreover, for any $ l>k, $ we also have
	\begin{align}\label{en54}
	\left| x \right|\leq k \ \text{and}  \ \left| y \right|\leq k\Rightarrow \left\{ \begin{array}{c}
	b_{i,k}\left( t,x,y \right) =b_{i,l}\left( t,x,y \right) =b_i\left( t,x,y \right) ,\\
	f_{1,k}\left( t,x \right) =f_{1,l}\left( t,x \right) =f_1\left( t,x \right), \qquad\ \\
	g_{1,k}\left( t,x \right) =g_{1,l}\left( t,x \right) =g_1\left( t,x \right),\qquad\ \\
	\end{array} \right. 
	\end{align}
	and
	\begin{align}\label{en55}
	\varLambda _{\alpha}^{0,t}\left( B^H \right) \leq k\Rightarrow  B^{H,k}_t= B^{H,l}_t= B_t^H.
	\end{align}
	Then, using the same argument as \cite{Guerra2008Stochastic}, it is easy to get that for any fixed $ k\in \mathbb{N}, $ there exists a unique strong solution $ (u_t^{\epsilon,k},v_t^{\epsilon,k} ) \in (W_0^{\alpha,\infty} (0,T;\mathbb{R}^n  ), W_0^{\alpha,\infty} (0,T;\mathbb{R}^m  ) )   $ to equation (\ref{orginal4}).

	\begin{lem}\label{lem3.3}
		Under the assumptions {\rm (A1)-(A4)}, for any $ \alpha \in \left(  1-H, \frac{1}{2} \land\beta  \land  \frac{\gamma}{2} \right), p\ge 1 $ and  fixed $ \epsilon_0\in \left( 0,1\right],   $ there exists some positive constant, such that   
		\begin{align}\label{en3.4}
		\mathbb{E}\left\|  u^{\epsilon_0 ,k} \right\|  _{\alpha ,\infty}^{p}\leq C_{\alpha,p,x, T}, \qquad 
		\end{align}
		and
		\begin{align}\label{en3.5}
		\mathbb{E}\left\|  v^{\epsilon_0 ,k}  \right\|  _{\alpha,\infty }^{p}\leq C_{\alpha,p,\epsilon_0,x,y, T}. 
		\end{align}
	\end{lem}
	\para{Proof:} First, we estimate $ \mathbb{E}\lVert  u^{\epsilon_0 ,k} \lVert _{\alpha ,\infty}^{p}. $ For brevity, we denote
	\begin{align}
	\varPsi _{t}^{\epsilon_0}\left( \lambda ,u^{\epsilon_0 ,k} \right) =\underset{r\in \left[ 0,t \right]}{\text{sup}}e^{-\lambda r}\left| u_{r}^{\epsilon_0 ,k} \right|,
	\end{align}
	and 
	\begin{align}
	\varPhi _{t}^{\epsilon_0}\left( \lambda ,u^{\epsilon_0 ,k} \right) =\underset{r\in \left[ 0,t \right]}{\text{sup}}e^{-\lambda r}\int_0^r{\frac{\left| u_{r}^{\epsilon_0 ,k}-u_{s}^{\epsilon_0 ,k} \right|}{\left( r-s \right) ^{\alpha +1}}ds}.
	\end{align}
	In order to estimate $ \lVert u^{\epsilon_0,k}\lVert_{\alpha,\infty},  $ we first estimate $ \varPsi _{t}^{\epsilon_0}\left( \lambda ,u^{\epsilon_0 ,k} \right). $ Thanks to the assumption (A4) and (\ref{en3}), it yields
	\begin{align}\label{en50}
	\varPsi _{t}^{\epsilon_0}\left( \lambda ,u^{\epsilon_0 ,k} \right) 
	%&=\underset{r\in \left[ 0,t \right]}{\text{sup}}e^{-\lambda r}\Big| x+\int_0^r{b_{1,k}\left( s,u_{s}^{\epsilon_0 ,k},v_{s}^{\epsilon_0 ,k} \right) ds}+\int_0^r{f_{1,k}\left( s,u_{s}^{\epsilon_0 ,k} \right) dW_{s}^{1}}  \cr
	%&\qquad\qquad\quad\quad\  +\int_0^r{g_{1,k}\left( s,u_{s}^{\epsilon_0 ,k} \right) dB_{s}^{H,k}} \Big|\cr
	&\le C_{x, T}\Big( 1+ \underset{r\in \left[ 0,t \right]}{\text{sup}}\Big| \int_0^r{f_{1,k}\left( s,u_{s}^{\epsilon_0 ,k} \right) dW_{s}^{1}} \Big|  \cr
	&\qquad\qquad   +\varLambda _{\alpha}^{0,t}\left( B^{H,k} \right) \underset{r\in \left[ 0,t \right]}{\text{sup}}\int_0^r{e^{-\lambda \left( r-s \right)}\left( s^{-\alpha}\varPsi _{s}^{\epsilon_0}\left( \lambda ,u^{\epsilon_0 ,k} \right) +\varPhi _{s}^{\epsilon_0}\left( \lambda ,u^{\epsilon_0 ,k} \right) \right)}ds \Big) \cr
	&\le C_{x, T}\left( 1+\varLambda _{\alpha}^{0,t}\left( B^{H,k} \right) \right) \Big( 1+ \lambda ^{\alpha -1}\varPsi _{t}^{\epsilon_0}\left( \lambda ,u^{\epsilon_0 ,k} \right) +\lambda ^{-1}\varPhi _{t}^{\epsilon_0}\left( \lambda ,u^{\epsilon_0 ,k} \right)  \cr
	& \qquad\qquad\quad\  +\underset{r\in \left[ 0,t \right]}{\text{sup}}\Big| \int_0^r{f_{1,k}\left( s,u_{s}^{\epsilon_0 ,k} \right) dW_{s}^{1}} \Big|  \Big),
	\end{align}
	where the last equation used the following estimate \cite[page 66]{Nualart2002Differential}
	\begin{align}
	\int_0^t{e^{-\lambda \left( t-r \right)}r^{-\alpha}}dr\leq C\lambda ^{\alpha -1}.\nonumber
	\end{align}
	Then, we estimate $ \varPhi _{t}^{\epsilon_0}\left( \lambda ,u^{\epsilon_0 ,k} \right). $ According to the equation (4.17) in \cite{Nualart2002Differential}, we can get
	\begin{align}\label{en51}
	\varPhi _{t}^{\epsilon_0}\left( \lambda ,u^{\epsilon_0 ,k} \right) 
	%&=\underset{r\in \left[ 0,t \right]}{\text{sup}}e^{-\lambda r}\Big| \int_0^r{\frac{\left| \int_s^r{b_{1,k}\left( \sigma ,u_{\sigma}^{\epsilon_0,k},v_{\sigma}^{\epsilon_0,k} \right)}d\sigma \right|}{\left( r-s \right) ^{\alpha +1}}}ds+\int_0^r{\frac{\left| \int_s^r{f_{1,k}\left( \sigma ,u_{\sigma}^{\epsilon_0,k} \right)}dW_{\sigma}^{1} \right|}{\left( r-s \right) ^{\alpha +1}}}ds 
	% \cr
	%&\qquad\qquad\qquad\    +\int_0^r{\frac{\left| \int_s^r{g_{1,k}\left( \sigma ,u_{\sigma}^{\epsilon_0,k} \right)}dB_{\sigma}^{H,k} \right|}{\left( r-s \right) ^{\alpha +1}}}ds \Big|
	%\cr
	&\le C_{\alpha ,x, T}\Big( 1+\underset{r\in \left[ 0,t \right]}{\text{sup}}\int_0^r{\frac{\left| \int_s^r{f_{1,k}\left( \sigma ,u_{\sigma}^{\epsilon_0,k} \right)}dW_{\sigma}^{1} \right|}{\left( r-s \right) ^{\alpha +1}}}ds+\varLambda _{\alpha}^{0,t}\left( B^{H,k} \right) \underset{r\in \left[ 0,t \right]}{\text{sup}}e^{-\lambda r} \cr
	  \times \Big(& \int_0^r{\frac{\left| g_{1,k}\left( s,u_{s}^{\epsilon_0,k} \right) \right|}{\left( r-s \right) ^{2\alpha}}ds} +\int_0^r{\int_0^s{\left( r-\sigma \right) ^{-\alpha}\frac{\left| g_{1,k}\left( s,u_{s}^{\epsilon_0,k} \right) -g_{1,k}\left( \sigma ,u_{\sigma}^{\epsilon_0,k} \right) \right|}{\left( s-\sigma \right) ^{1+\alpha}}d\sigma ds}} \Big) \Big) \cr
	&\le C_{\alpha ,x, T}\left( 1+\varLambda _{\alpha}^{0,t}\left( B^{H,k} \right) \right) \Big( 1+\underset{r\in \left[ 0,t \right]}{\text{sup}}\int_0^r{\left( r-s \right) ^{-\alpha -1}\Big| \int_s^r{f_{1,k}\left( \sigma ,u_{\sigma}^{\epsilon_0,k} \right)}dW_{\sigma}^{1} \Big|}ds   \cr
	& \qquad\qquad\qquad\qquad\qquad\quad\quad\ \  +\lambda ^{2\alpha -1}\varPsi _{t}^{\epsilon_0}\left( \lambda ,u^{\epsilon_0 ,k} \right) +\lambda ^{\alpha -1}\varPhi _{t}^{\epsilon_0}\left( \lambda ,u^{\epsilon_0 ,k} \right) \Big), 
	\end{align}
	where the last equation used the estimate \cite[page 66]{Nualart2002Differential}
	\begin{align}
	\int_0^r{e^{-\lambda \left( r-s \right)}\left( r-s \right) ^{-\alpha}}ds
	%=\frac{1}{\lambda}\int_0^{\lambda r}{e^{-s}\lambda ^{\alpha}s^{-\alpha}}ds\le \lambda ^{\alpha -1}\int_0^{\infty}{e^{-s}s^{-\alpha}}ds
	=\varGamma \left( 1-\alpha \right) \lambda ^{\alpha -1}.
	\nonumber
	\end{align}
	
	Let $ C:=C_{x, T}\vee C_{\alpha ,x, T} $ and $ \lambda=\left( 4C\left( 1+\varLambda _{\alpha}^{0,t}\left( B^{H,k} \right) \right)\right)^{\frac{1}{1-\alpha}}.   $
	Combining (\ref{en50}) and (\ref{en51}), making simple transformations as  \cite{pei2020averaging},   it is easy to get that  
	\begin{align}
	\varPsi _{t}^{\epsilon_0}\left( \lambda ,u^{\epsilon_0 ,k} \right) +\varPhi _{t}^{\epsilon_0}\left( \lambda ,u^{\epsilon_0 ,k} \right) &\le C_{\alpha ,x, T}\left( 1+\varLambda _{\alpha}^{0,t}\left( B^{H,k}\right) \right) ^{\frac{1}{1-\alpha}}\left(  1+\varUpsilon _{t}^{\epsilon_0}\left( \lambda ,u^{\epsilon_0 ,k} \right)\right),\nonumber
	\end{align}
	where
	\begin{align}
	\varUpsilon _{t}^{\epsilon_0}\left( \lambda ,u^{\epsilon_0 ,k} \right) :=\underset{r\in \left[ 0,t \right]}{\text{sup}}\Big( \Big| \int_0^r{f_{1,k}\left( s,u_{s}^{\epsilon_0 ,k} \right) dW_{s}^{1}} \Big|+ \int_0^r{\left( r-s \right) ^{-\alpha -1}\Big| \int_s^r{f_{1,k}\left( \sigma ,u_{\sigma}^{\epsilon_0,k} \right)}dW_{\sigma}^{1} \Big|}ds\Big).\nonumber
	\end{align}
	Hence
	\begin{align}
	\left\|  u^{\epsilon_0 ,k} \right\|  _{\alpha ,\infty}&\le e^{\lambda T}\left( \varPsi _{T}^{\epsilon_0}\left( \lambda ,u^{\epsilon_0 ,k} \right) +\varPhi _{T}^{\epsilon_0}\left( \lambda ,u^{\epsilon_0 ,k} \right) \right) \cr
	%&\le C_{\alpha ,x, T}e^{\left( 4C\left( 1+\varLambda _{\alpha}^{0,t}\left( B^H \right) \right) \right) ^{\frac{1}{1-\alpha}}T}\left( 1+\varLambda _{\alpha}^{0,t}\left( B^H \right) \right) ^{\frac{1}{1-\alpha}}\left( 1+\varUpsilon _{t}^{\epsilon_0}\left( \lambda ,u^{\epsilon_0 ,k} \right) \right) \cr
	&\le C_{\alpha ,x, T}e^{C_{\alpha ,x, T}\left( \varLambda _{\alpha}^{0,T}\left( B^{H,k} \right) \right) ^{\frac{1}{1-\alpha}}}\left( 1+\varLambda _{\alpha}^{0,T}\left( B^{H,k}  \right) \right) ^{\frac{1}{1-\alpha}}\left( 1+\varUpsilon _{T}^{\epsilon_0}\left( \lambda ,u^{\epsilon_0 ,k} \right) \right).\nonumber
	\end{align}
	By taking expectations on both sides of the above equation, we can get
	\begin{align}\label{en52}
	\mathbb{E}\left\|  u^{\epsilon_0 ,k} \right\|  _{\alpha ,\infty}^{p}&\le C_{\alpha ,x, T}\Big[ \mathbb{E}\big[ e^{3pC_{\alpha ,x, T}\left( \varLambda _{\alpha}^{0,T}\left( B^{H,k} \right) \right) ^{\frac{1}{1-\alpha}}} \big] \big( 1+\mathbb{E}\left[  \varLambda _{\alpha}^{0,T}\left( B^{H,k} \right) \right]  ^{\frac{3p}{1-\alpha}} \big)\cr
	&\qquad\qquad\qquad\times \left( 1+\mathbb{E}\left[  \varUpsilon _{T}^{\epsilon_0}\left( \lambda ,u^{\epsilon_0 ,k} \right) \right]  ^{3p} \right) \Big] ^{\frac{1}{3}}.
	\end{align}
	We need to prove that $ \mathbb{E}\left[  \varUpsilon _{T}^{\epsilon_0}\left( \lambda ,u^{\epsilon_0 ,k} \right) \right]  ^{p}  $ is bounded for any $ p\geq 1. $ Applying the Garsia-Rodemich-Rumsey inequality (\ref{en41}) with $ p=2/\eta  $ and $ \theta=\left( 1-\eta \right)/2  $ (where $ \eta\in \left(0, 1/2-\alpha \right) $),  it deduce that
	\begin{align}\label{en57}
	\Big| \int_s^t{f_{1,k}\left( r,u_{r}^{\epsilon_0,k} \right)}dW_{r}^{1} \Big|\le C_{\eta  }\left| t-s \right|^{1/2-\eta}\zeta, 
	\end{align}
	where
	\begin{align}
	\zeta=\Big( \int_0^T{\int_0^T{\left| \sigma -r \right|^{-1/\eta}\Big| \int_r^{\sigma}{f_{1,k}\left( z,u_{z}^{\epsilon_0,k} \right)}dW_{z}^{1} \Big|^{2/\eta}d\sigma dr}} \Big) ^{\eta /2}.\nonumber
	\end{align}
	Then, taking expectations for $ \zeta, $ we can get
	\begin{align}
	\mathbb{E}\zeta ^p
	%&=\mathbb{E}\Big( \int_0^T{\int_0^T{\left| \sigma -r \right|^{-1/\eta}\Big| \int_r^{\sigma}{f_1\left( z,u_{z}^{\epsilon_0,k} \right)}dW_{z}^{1} \Big|^{2/\eta}d\sigma dr}} \Big) ^{\left( \eta p \right) /2}\cr
	%&=\int_0^T{\int_0^T{\left| \sigma -r \right|^{-p/2}\mathbb{E}\Big| \int_r^{\sigma}{f_{1,k}\left( z,u_{z}^{\epsilon_0,k} \right)}dW_{z}^{1} \Big|^pd\sigma dr}}\cr
	= \int_0^T{\int_0^T{\left| \sigma -r \right|^{-p/2}\mathbb{E}\Big( \int_r^{\sigma}{\left| f_{1,k}\left( z,u_{z}^{\epsilon_0,k} \right) \right|^2}dz \Big) ^{p/2}d\sigma dr}} 
	%&\le C_{p,\eta ,T}\int_0^T{\int_0^T{\left| \sigma -r \right|^{-p/2}\Big( \int_r^{\sigma}{\mathbb{E}\left| f_1\left( z,u_{z}^{\epsilon_0} \right) -f_1\left( t_2,u_{z}^{\epsilon_0} \right) \right|^2+\mathbb{E}\left| f_1\left( t_2,u_{z}^{\epsilon_0} \right) \right|^2}dz \Big) ^{p/2}d\sigma dr}}\cr
	\le C_{T}.\nonumber
	\end{align}
	Fixed $ \eta\in \left(0, 1/2 -\alpha\right),  $ it follows that
	\begin{eqnarray}\label{en53}
	\mathbb{E}\left[  \varUpsilon _{T}^{\epsilon_0}\left( \lambda ,u^{\epsilon_0 ,k} \right) \right]  ^{p}
	%&=\mathbb{E}\underset{t\in \left[ 0,T \right]}{\text{sup}}\left\|  \int_0^t{f_1\left( r,u_{r}^{\epsilon_0} \right)}dW_{r}^{1} \right\|  _{\alpha}^{p}\cr
	&\leq& C_p\int_0^T{\mathbb{E}\left| f_{1,k}\left( r,u_{r}^{\epsilon_0,k} \right)   \right|^p }dr +C_p\mathbb{E}\zeta ^p\Big( \underset{t\in \left[ 0,T \right]}{\text{sup}}\Big| \int_0^t{\left( t-s \right) ^{-\eta -\alpha -1/2}}ds \Big|^p \Big) \cr 
	&\le& C_{\alpha,p, T}. 
	\end{eqnarray}
	Thanks to (\ref{en52}), (\ref{en53}) and $ \Lambda_{\alpha}^{0,T}\left(B^{H,k}\right) 
	$ has moments of all orders \cite{Nualart2002Differential}, it yields (\ref{en3.4}).
	
	Next, we estimate $ \mathbb{E}\lVert   v^  {\epsilon_0 ,k}_t \lVert _{\alpha,\infty }^{2}   $ 
	\begin{align}\label{en60}
	\mathbb{E}\left\|   v^  {\epsilon_0 ,k} \right\|  _{\alpha ,\infty}^{p}
	%&=\mathbb{E}\Big( \underset{t\in \left[ 0,T \right]}{\text{sup}}\Big\|  y+\frac{1}{\epsilon_0}\int_0^t{b_{2,k}\left( r,u_{r}^{\epsilon_0 ,k},v_{r}^{\epsilon_0 ,k} \right) dr}+\frac{1}{\sqrt{\epsilon_0}}\int_0^t{f_2\left( r,v_{r}^{\epsilon_0 ,k} \right) dW_{r}^{2}} \Big\| _{\alpha} \Big) ^p\cr
	%&\le C_{p,y}\Big( 1+\mathbb{E}\Big[ \underset{t\in \left[ 0,T \right]}{\text{sup}}\Big\|  \int_0^t{b_{2,k}\left( r,u_{r}^{\epsilon_0 ,k},v_{r}^{\epsilon_0 ,k} \right) dr} \Big\| _{\alpha}\Big] ^{p}+\mathbb{E}\Big[ \underset{t\in \left[ 0,T \right]}{\text{sup}}\Big\|  \int_0^t{f_2\left( r,v_{r}^{\epsilon_0 ,k} \right) dW_{r}^{1}} \Big\|  _{\alpha}\Big] ^{p} \Big) \cr
	&\le C_{p,y,\epsilon_0}\Big( 1+\mathbb{E}\Big[ \underset{t\in \left[ 0,T \right]}{\text{sup}}\Big| \int_0^t{b_{2,k}\left( r,u_{r}^{\epsilon_0 ,k},v_{r}^{\epsilon_0 ,k} \right) dr} \Big|\Big] ^p+  \varGamma _{T}^{\epsilon_0}\left( b_{2,k} \right) \cr
	&  \qquad\qquad+\mathbb{E}\Big[ \underset{t\in \left[ 0,T \right]}{\text{sup}}\Big| \int_0^t{f_2\left( r,v_{r}^{\epsilon_0 ,k} \right) dW_{r}^{2}} \Big|\Big] ^p+\varLambda _{T}^{\epsilon_0}\left( f_2 \right) \Big) \cr
	&\le C_{p,y,\epsilon_0,T}\Big( 1+\int_0^T{\big(\mathbb{E}\left| u_{r}^{\epsilon_0 ,k} \right|^p+\mathbb{E}\left| v_{r}^{\epsilon_0 ,k} \right|^p\big)dr} +\varGamma _{T}^{\epsilon_0}\left( b_{2,k} \right) +\varLambda _{T}^{\epsilon_0}\left( f_2 \right) \Big),
	\end{align}
	where
	\begin{align}
	\varGamma _{T}^{\epsilon_0}\left( b_{2,k}\right) :=\mathbb{E}\Big( \underset{t\in \left[ 0,T \right]}{\text{sup}}\int_0^t{{\left( t-s \right) ^{-\alpha -1}{\Big| \int_s^t{b_{2,k}\left( r,u_{r}^{\epsilon_0 ,k},v_{r}^{\epsilon_0 ,k} \right)}dr \Big|}}}ds \Big) ^p,\nonumber
	\end{align}
	and 
	\begin{align}
	\varLambda _{T}^{\epsilon_0}\left( f_2 \right) :=\mathbb{E}\Big(\underset{t\in \left[ 0,T \right]}{\text{sup}} \int_0^t{{\left( t-s \right) ^{-\alpha -1}}{\Big| \int_s^t{f_2\left( r,v_{r}^{\epsilon_0 ,k} \right)}dW_{r}^{2} \Big|}}ds \Big) ^p.\nonumber
	\end{align}
	For $ \varGamma _{T}^{\epsilon_0}\left( b_{2,k} \right), $ applying the Garsia-Rodemich-Rumsey inequality (\ref{en41}) with $ p=2/\varrho  $ and $ \theta=\left( 1-\varrho \right)/2  $ (where $ \varrho\in \left(0, 1/2-\alpha \right) $) again,  we also can get
	\begin{align}
	\Big| \int_s^t{b_{2,k}\left( r,u_{r}^{\epsilon_0 ,k},v_{r}^{\epsilon_0 ,k} \right)}dr \Big|\le C_{\varrho  }\left| t-s \right|^{1/2-\varrho}\varsigma,\nonumber
	\end{align}
	where
	\begin{align}
	\varsigma=\Big( \int_0^T{\int_0^T{\left| \sigma -r \right|^{-1/\varrho}\Big| \int_r^{\sigma}{b_{2,k}\left( z,u_{z}^{\epsilon_0 ,k},v_{z}^{\epsilon_0 ,k} \right)}dz \Big|^{2/\varrho}d\sigma dr}} \Big) ^{\varrho /2}.\nonumber
	\end{align}
	Then, taking expectations for $ \varsigma,  $ we can get
	\begin{align}
	\mathbb{E}\varsigma ^p
	%&=\mathbb{E}\Big( \int_0^T{\int_0^T{\left| \sigma -r \right|^{-1/\eta}\Big| \int_r^{\sigma}{f_1\left( z,u_{z}^{\epsilon_0,k} \right)}dW_{z}^{1} \Big|^{2/\eta}d\sigma dr}} \Big) ^{\left( \eta p \right) /2}\cr
	%&=\int_0^T{\int_0^T{\left| \sigma -r \right|^{-p/2}\mathbb{E}\Big| \int_r^{\sigma}{b_{2,k}\left( z,u_{z}^{\epsilon_0 ,k},v_{z}^{\epsilon_0 ,k} \right)}dz \Big|^pd\sigma dr}}\cr
	&\le \int_0^T{\int_0^T{\left| \sigma -r \right|^{-1+p/2}  \int_r^{\sigma}{\mathbb{E}\left|b_{2,k}\left( z,u_{z}^{\epsilon_0 ,k},v_{z}^{\epsilon_0 ,k} \right)  \right|^p}dz d\sigma dr}}\cr
	&\le \int_0^T{\int_0^T{\left| \sigma -r \right|^{-1+p/2} d\sigma dr \int_0^{T}{\mathbb{E}\left|b_{2,k}\left( z,u_{z}^{\epsilon_0 ,k},v_{z}^{\epsilon_0 ,k} \right)  \right|^p}dz }}\cr
	%&=\int_0^T{\int_0^z{\int_z^T{\left| \sigma -r \right|^{-1+p/2}d\sigma dr}\mathbb{E}\left| b_{2,k}\left( z,u_{z}^{\epsilon_0 ,k},v_{z}^{\epsilon_0 ,k} \right) \right|^pdz}}\cr
	&\le C_{p,T}\Big( 1+\int_0^T{\mathbb{E}\left| u_{r}^{\epsilon_0 ,k} \right|^p}dr+\int_0^T{\mathbb{E}\left| v_{r}^{\epsilon_0 ,k} \right|^p}dr
	\Big) .\nonumber
	\end{align}
	Hence 
	\begin{align}\label{en58}
	\varGamma _{T}^{\epsilon_0}\left( b_{2,k} \right)&\le C\mathbb{E}\varsigma^p\Big( \underset{t\in \left[ 0,T \right]}{\text{sup}}\int_0^t{\left( t-s \right) ^{-\varrho -\alpha -1/2}}ds \Big) ^p\cr
	&\le  C_{\alpha,p,T}\Big(1+\int_0^T{\mathbb{E}\left| u_{r}^{\epsilon_0 ,k} \right|^p}dr+\int_0^T{\mathbb{E}\left| v_{r}^{\epsilon_0 ,k} \right|^p}dr\Big).
	\end{align}
	Moreover, use the same argument as (\ref{en53}) and (\ref{en58}), we also can obtain
	\begin{align}\label{en59} 
	\varLambda_{T}^{\epsilon_0}\left( f_2 \right)
	&\le  C_{\alpha,p,T}\Big(1+\int_0^T{\mathbb{E}\left| v_{r}^{\epsilon_0 ,k} \right|^p}dr\Big).
	\end{align}
	Substituting (\ref{en58}) and (\ref{en59}) into (\ref{en60}), thanks to (\ref{en3.4}), we have
	\begin{align}
	\mathbb{E}\lVert  v^  {\epsilon_0 ,k} \lVert  _{\alpha ,\infty}^{p}=\mathbb{E}\Big( \underset{t\in \left[ 0,T \right]}{\text{sup}}\lVert v_{t}^{\epsilon_0 ,k} \lVert  _{\alpha} \Big) ^p\le C_{\alpha ,p,\epsilon_0,x,y,T}\Big( 1+\int_0^T{\mathbb{E}\Big( \underset{\sigma \in \left[ 0,r \right]}{\text{sup}}\lVert v_{\sigma}^{\epsilon_0 ,k} \lVert  _{\alpha} \Big) ^pdr} \Big).
	\end{align}
	Then, by Gronwall inequality, we have (\ref{en3.5}). The proof is completed.\qed

	\subsection{The existence and uniqueness of solutions}
	Now, we study the existence and uniqueness of solutions for original equation (\ref{orginal1}):
	\para{Proof of \thmref{thm3.1}:} In order to prove the existence of the solution for (\ref{orginal1}), fixed $ \epsilon_0\in \left( 0,1\right]    $ and any $ k\in \mathbb{N}, $ we define the following stopping time 
	\begin{align}
	\tau _k^1:=\text{inf}\left\{ t\geq 0:\lVert  u_{t}^{\epsilon_0 ,k} \lVert  _{\alpha}+\lVert v_{t}^{\epsilon_0 ,k} \lVert  _{\alpha} \geq k \right\} \land \tau_k,
	\end{align}
	and we let
	\begin{align}
	\tau :=\underset{k\in \mathbb{N}}{\text{sup}}\tau _k^1.
	\end{align}
	It is easy to know that the sequence of stopping times $ \left\lbrace \tau _k^1\right\rbrace  $ is non-decreasing and $ \mathbb{P}(\tau =+\infty)=1 $. Indeed,
	$$
	\mathbb{P}\left( \tau <+\infty \right) =\lim_{T\rightarrow +\infty}\mathbb{P}\left( \tau \leq T \right), 
	$$
	and for each $ T > 0, $ thanks to \lemref{lem3.3} and Lemma 7.5 in \cite{Nualart2002Differential},  we can get
	\begin{align}
	\mathbb{P}\left( \tau  \leq T \right) &=\lim_{k\rightarrow +\infty}\mathbb{P}\left( \tau _{k}^1 \leq T \right)\cr
	&=\lim_{k\rightarrow +\infty}\mathbb{P}\Big(\sup_{t\in \left[ 0,T \right]}\lVert  u_{t}^{\epsilon_0 ,k} \lVert  _{\alpha}+\sup_{t\in \left[ 0,T \right]}\lVert  v_{t}^{\epsilon_0 ,k} \lVert  _{\alpha}+ \sup_{t\in \left[ 0,T \right]} \varLambda _{\alpha}^{0,t}\left( B^H  \right) \geq k \Big)\cr
	&=\lim_{k\rightarrow +\infty}\frac{1}{k^2}\mathbb{E}\Big(\sup_{t\in \left[ 0,T \right]}\lVert  u_{t}^{\epsilon_0 ,k} \lVert _{\alpha}^2+\sup_{t\in \left[ 0,T \right]}\lVert  v_{t}^{\epsilon_0 ,k} \lVert  _{\alpha}^2+\frac{1}{\Gamma^2(1-\alpha)} \sup _{0<s<t<T}\left|\left(D_{t-}^{1-\alpha} B^H _{t-}\right)(s)\right|^2  \Big) \cr
	&=0. \nonumber
	\end{align}
	Hence, $ P(\tau <+\infty)=0, $ that is, $ P(\tau =+\infty)=1. $ Further, for any $ t\in [0,T] $ and $ \omega\in \left\lbrace \tau   =+\infty\right\rbrace  $, there exists $ l\in \mathbb{N} $ such that $ t\leq \tau _{l}^1 (\omega)$. Then, we define
	\begin{align}
	u_t^{\epsilon_0}\left(\omega\right):= u_t^{\epsilon_0,l}\left(\omega\right)\quad \text{and} \quad v_t^{\epsilon_0}\left(\omega\right):= v_t^{\epsilon_0,l}\left(\omega\right).
	\end{align}
	This is a good definition, as for any $  t\leq \tau _{k}^1 \land\tau _{l}^1,   $
	we have 
	\begin{align}\label{en56}
	u_t^{\epsilon_0,k}  = u_t^{\epsilon_0,l}\quad  \text{and} \quad v_t^{\epsilon_0,k}  = v_t^{\epsilon_0,l}, \quad \mathbb{P} -a.s. 
	\end{align}  
	Actually, for any  
	$ k\geq l  $ and $ t\leq \tau _{k}^1 \land\tau _{l}^1, $ thanks to (\ref{en54}) and (\ref{en55}),  we have
	\begin{align}
	u_{t}^{\epsilon_0 ,k}-u_{t}^{\epsilon_0 ,l}
	&=\int_0^t{\left[ b_{1,k}\left( r,u_{r}^{\epsilon_0 ,k},v_{r}^{\epsilon_0 ,k} \right) -b_{1,l}\left( r,u_{r}^{\epsilon_0 ,l},v_{r}^{\epsilon_0 ,l} \right) \right]}dr \cr
	&\quad+\int_0^t{\left[ f_{1,k}\left( r,u_{r}^{\epsilon_0 ,k} \right) -f_{1,l}\left( r,u_{r}^{\epsilon_0 ,l} \right) \right]}dW_r^1\cr
	&\quad+\int_0^t{g_{1,k}\left( r,u_{r}^{\epsilon_0 ,k} \right)}dB_{r}^{H,k}-\int_0^t{g_{1,l}\left( r,u_{r}^{\epsilon_0 ,l} \right)}dB_{r}^{H,l}\cr 
	&=\int_0^t{\left[ b_{1,k}\left( r,u_{r}^{\epsilon_0 ,k},v_{r}^{\epsilon_0 ,k} \right) -b_{1,k}\left( r,u_{r}^{\epsilon_0 ,l},v_{r}^{\epsilon_0 ,l} \right) \right]}dr \cr
	&\quad+\int_0^t{\left[ f_{1,k}\left( r,u_{r}^{\epsilon_0 ,k} \right) -f_{1,k}\left( r,u_{r}^{\epsilon_0 ,l} \right) \right]}dW_r^1\cr
	&\quad+\int_0^t{\left[ g_{1,k}\left( r,u_{r}^{\epsilon_0 ,k} \right) -g_{1,k}\left( r,u_{r}^{\epsilon_0 ,l} \right) \right]}dB_{r}^{H,k}.\nonumber
	\end{align}
	Similarly, we also have
	\begin{align}
	v_{t}^{\epsilon_0 ,k}-v_{t}^{\epsilon_0 ,l}
	%&=\frac{1}{\epsilon_0}\int_0^t{\left[ b_{2,k}\left( r,u_{r}^{\epsilon_0 ,k},v_{r}^{\epsilon_0 ,k} \right) -b_{2,l}\left( r,u_{r}^{\epsilon_0 ,l},v_{r}^{\epsilon_0 ,l} \right) \right]}dr\cr
	%&\quad+\frac{1}{\sqrt{\epsilon_0}}\int_0^t{\left[ f_2\left( r,u_{r}^{\epsilon_0 ,k},v_{r}^{\epsilon_0 ,k} \right) -f_2\left( r,u_{r}^{\epsilon_0 ,l},v_{r}^{\epsilon_0 ,l} \right) \right]}dW_r^2\cr 
	&=\frac{1}{\epsilon_0}\int_0^t{\left[ b_{2,k}\left( r,u_{r}^{\epsilon_0 ,k},v_{r}^{\epsilon_0 ,k} \right) -b_{2,k}\left( r,u_{r}^{\epsilon_0 ,l},v_{r}^{\epsilon_0 ,l} \right) \right]}dr\cr
	&\quad+\frac{1}{\sqrt{\epsilon_0}}\int_0^t{\left[ f_2\left( r,u_{r}^{\epsilon_0 ,k},v_{r}^{\epsilon_0 ,k} \right) -f_2\left( r,u_{r}^{\epsilon_0 ,l},v_{r}^{\epsilon_0 ,l} \right) \right]}dW_r^2.\nonumber
	\end{align}
	According to the paper \cite{Nualart2002Differential} and \cite{Guerra2008Stochastic}, we know that the trajectories of  $ u_{t}^{\epsilon_0 ,k} $ and $ u_{t}^{\epsilon_0 ,l} $ are $ \eta $-H\"{o}lder continuous for all $ \eta< 1/2. $ Now, let $ \eta\in  ( \frac{\alpha}{\gamma}, \frac{1}{2})  $ and consider the set $ \Omega_N\subset\Omega $ with $ N\in\mathbb{N}, $ such that
	\begin{align}
	\Omega_N:=\big\lbrace \omega\in \left\lbrace \tau   =+\infty\right\rbrace:\left\| u^{\epsilon_0 ,k}\right\| _\eta  \leq N \ \text{and}\ \left\| u^{\epsilon_0 ,l}\right\| _\eta \leq N \big\rbrace. \nonumber
	\end{align}
	It is clear that $  \Omega_N\nearrow  \left\lbrace \tau   =+\infty\right\rbrace. $ Then, by proceeding as Proposition 3.4, Proposition 3.6 and Proposition 3.9 in \cite{Guerra2008Stochastic}, we can get 
	\begin{align}\label{en525}
	&\qquad\mathbb{E}\big[ \big\|  u_{t}^{\epsilon_0 ,k}-u_{t}^{\epsilon_0 ,l} \big\|  _{\alpha}^{2}\mathbf{1}_{\varOmega _N} \big] +\mathbb{E}\big[ \big\|  v_{t}^{\epsilon_0 ,k}-v_{t}^{\epsilon_0 ,l} \big\|  _{\alpha}^{2}\mathbf{1}_{\varOmega _N} \big] \cr
	&\le C_{\alpha ,T}\mathbb{E}\Big( \int_0^t{\left( t-r \right) ^{-\alpha}\left\|  b_{1,k}\left( r,u_{r}^{\epsilon_0 ,k},v_{r}^{\epsilon_0 ,k} \right) -b_{1,k}\left( r,u_{r}^{\epsilon_0 ,l},v_{r}^{\epsilon_0 ,l} \right) \right\|  _{\alpha}\mathbf{1}_{\varOmega _N}}dr \Big) ^2\cr
	&\quad+C_{\alpha ,T}\int_0^t{\left( t-r \right) ^{-\frac{1}{2}-\alpha}\mathbb{E}\big[ \left\|  f_{1,k}\left( r,u_{r}^{\epsilon_0 ,k} \right) -f_{1,k}\left( r,u_{r}^{\epsilon_0 ,l} \right) \right\|  _{\alpha}^{2}\mathbf{1}_{\varOmega _N} \big]}dr\cr
	&\quad+C_{\alpha ,T}\mathbb{E}\Big[ \varLambda _{\alpha}^{0,t}\left( B^{H,k} \right) \int_0^t{\left( \left( t-r \right) ^{-2\alpha}+r^{-\alpha} \right) \left( 1+\Delta u_{r}^{\epsilon_0 ,k}+\Delta u_{r}^{\epsilon_0 ,l} \right) \left\|  u_{r}^{\epsilon_0 ,k}-u_{r}^{\epsilon_0 ,l} \right\|  _{\alpha}\mathbf{1}_{\varOmega _N}}dr \Big] ^2\cr
	&\quad+C_{\alpha,\epsilon_0 ,T}\mathbb{E}\Big( \int_0^t{\left( t-r \right) ^{-\alpha}\left\|  b_{2,k}\left( r,u_{r}^{\epsilon_0 ,k},v_{r}^{\epsilon_0 ,k} \right) -b_{2,k}\left( r,u_{r}^{\epsilon_0 ,l},v_{r}^{\epsilon_0 ,l} \right) \right\|  _{\alpha}\mathbf{1}_{\varOmega _N}}dr \Big) ^2\cr
	&\quad+C_{\alpha,\epsilon_0 ,T}\int_0^t{\left( t-r \right) ^{-\frac{1}{2}-\alpha}\mathbb{E}\big[ \left\|  f_2\left( r,u_{r}^{\epsilon_0 ,k},v_{r}^{\epsilon_0 ,k} \right) -f_2\left( r,u_{r}^{\epsilon_0 ,l},v_{r}^{\epsilon_0 ,l} \right) \right\|  _{\alpha}^{2}\mathbf{1}_{\varOmega _N} \big]}dr\cr
	&\le C_{\alpha,\epsilon_0 ,T}\int_0^t{\big( \left( t-r \right) ^{-2\alpha}+\left( t-r \right) ^{-\frac{1}{2}-\alpha}\big) \left( \mathbb{E}\big[ \left\|  u_{r}^{\epsilon_0 ,k}-u_{r}^{\epsilon_0 ,l} \right\|  _{\alpha}^{2}\mathbf{1}_{\varOmega _N} \big] +\mathbb{E}\big[ \left\|  v_{r}^{\epsilon_0 ,k}-v_{r}^{\epsilon_0 ,l} \right\|  _{\alpha}^{2}\mathbf{1}_{\varOmega _N} \big] \right)}dr\cr
	%&\quad+C_{\alpha,\epsilon_0 ,T}\int_0^t{\left( t-r \right) ^{-\frac{1}{2}-\alpha}\left( \mathbb{E}\left[ \left\|  u_{r}^{\epsilon_0 ,k}-u_{r}^{\epsilon_0 ,l} \right\|  _{\alpha}^{2}\mathbf{1}_{\varOmega _N} \right] +\mathbb{E}\left[ \left\|  v_{r}^{\epsilon_0 ,k}-v_{r}^{\epsilon_0 ,l} \right\|  _{\alpha}^{2}\mathbf{1}_{\varOmega _N} \right] \right)}dr\cr
	&\quad+C_{\alpha ,k,k,T}\int_0^t{\left( \left( t-r \right) ^{-2\alpha}+r^{-\alpha} \right) \mathbb{E}\big[ \left\|  u_{r}^{\epsilon_0 ,k}-u_{r}^{\epsilon_0 ,l} \right\|  _{\alpha}^{2}\mathbf{1}_{\varOmega _N} \big]}dr\cr
	&\le C_{\alpha ,\epsilon_0,k,k,T}\int_0^t{\left( t-r \right) ^{-\frac{1}{2}-\alpha}r^{-\frac{1}{2}-\alpha}\left( \mathbb{E}\big[ \left\|  u_{r}^{\epsilon_0 ,k}-u_{r}^{\epsilon_0 ,l} \right\|  _{\alpha}^{2}\mathbf{1}_{\varOmega _N} \big] +\mathbb{E}\big[ \left\|  v_{r}^{\epsilon_0 ,k}-v_{r}^{\epsilon_0 ,l} \right\|  _{\alpha}^{2}\mathbf{1}_{\varOmega _N} \big] \right)}dr, \cr
	&\ 
	\end{align}
	where the last estimate is because $ \varLambda _{\alpha}^{0,t}\left( B^{H,k} \right)\leq k $ and  
	\begin{align}
	1+\Delta u_{r }^{\epsilon_0 ,k}+\Delta u_{r }^{\epsilon_0 ,l}&=1+\int_0^r{\frac{\left| u_{r }^{\epsilon_0 ,k}-u_{s }^{\epsilon_0 ,k} \right|^{\gamma}}{\left( r-s \right) ^{\alpha +1}}ds}+\int_0^r{\frac{\left| u_{r }^{\epsilon_0 ,l}-u_{s }^{\epsilon_0 ,l} \right|^{\gamma}}{\left( r-s \right) ^{\alpha +1}}ds}\cr
	&\le 1+\left( \left\|  u^{\epsilon_0 ,k} \right\|  _{\eta}^{\gamma}+\left\|  u^{\epsilon_0 ,l} \right\|  _{\eta}^{\gamma} \right) \int_0^r{\left( r-s \right) ^{\eta \gamma -\alpha -1}ds}\le C_N,\nonumber
	\end{align}
	as $ t\leq \tau _{k}^1\land\tau _{l}^1 $ and $ \Omega_N\subset  \left\lbrace \tau   =+\infty\right\rbrace. $ 
	Therefore, for any $ k\in \mathbb{N}, $ by the Gronwall-type lemma (Lemma 7.6 in  \cite{Nualart2002Differential}), we deduce that
	\begin{align}
	\mathbb{E}\big[ \lVert u_{t}^{\epsilon_0 ,k}-u_{t}^{\epsilon_0 ,l} \lVert  _{\alpha}^{2}\mathbf{1}_{\varOmega _N} \big] +\mathbb{E}\big[ \lVert  v_{t}^{\epsilon_0 ,k}-v_{t}^{\epsilon_0 ,l} \lVert  _{\alpha}^{2}\mathbf{1}_{\varOmega _N} \big] 
	=0, \quad  t\leq \tau _{k}\land\tau _{l}.
	\end{align}
	Then, let $ N\rightarrow+\infty, $ as $  \Omega_N\nearrow  \left\lbrace \tau   =+\infty\right\rbrace, $ we can get (\ref{en56}).
	
	Recalling that if $ \omega\in\left\lbrace \tau   =+\infty\right\rbrace $ and $ t \leq \tau _{l}, $ we denote $ u_t^{\epsilon_0} $ is equal to $ u_t^{\epsilon_0,l}  $ and $ v_t^{\epsilon_0} $ is equal to $ v_t^{\epsilon_0,l}, $ thanks to (\ref{en54}) and (\ref{en55}), it follows that
	\begin{align}
	\left\{ \begin{array}{l}
	u_{t}^{\epsilon_0}=x+\int_0^t{b_1\left( r,u_{r}^{\epsilon_0},v_{r}^{\epsilon_0} \right)}dr+\int_0^t{f_1\left( r,u_{r}^{\epsilon_0} \right)}dW_{r}^{1}+\int_0^t{g_1\left( r,u_{r}^{\epsilon_0} \right)}dB_{r}^{H},\\
	v_{t}^{\epsilon_0}=y+\frac{1}{\epsilon_0}\int_0^t{b_2\left( r,u_{r}^{\epsilon_0},v_{r}^{\epsilon_0} \right)}dr+\frac{1}{\sqrt{\epsilon_0}}\int_0^t{f_2\left( r,u_{r}^{\epsilon_0},v_{r}^{\epsilon_0} \right)}dW_{r}^{2},\\
	\end{array} \right.  
	\end{align}
	$ \mathbb{P} $-a.s., that is, for any fixed $  \epsilon_0\in (0,1], \ \left( u_{t}^{\epsilon_0}, v_{t}^{\epsilon_0}\right)  $ is a solution of equation (\ref{orginal1}).
	
	Finally, denote another solution of equation (\ref{orginal1}) is $ \left( u_{t}^{\epsilon_0,*}, v_{t}^{\epsilon_0,*}\right),  $ with the same argument as (\ref{en56}), we can also get that 
	\begin{align} 
	u_t^{\epsilon_0 }  = u_t^{\epsilon_0,*}\quad  \text{and} \quad v_t^{\epsilon_0 }  = v_t^{\epsilon_0,*}, \quad \mathbb{P} -a.s. 
	\end{align}
	Thus, we prove the uniqueness of the solution of equation (\ref{orginal1}). This proof is completed.\qed

	\section{Proof of  \thmref{thm2.1}}\label{sec-4}
	In this section, we prove the main \thmref{thm2.1}, i.e. the slow process $ u^\epsilon_t  $ converges to the averaged process $ \bar{u}_t  $ in the mean square sense, as $ \epsilon \rightarrow 0.$  
	Firstly, we need to define the averaged equation and give some properties of the averaged coefficient. Secondly, we construct an auxiliary process $ \hat{v}  ^\epsilon_t $  by the technique of time discretization and give some estimates about it on the basis of some a-priori estimates for the solution $ \left(  u^\epsilon_t,  v^  \epsilon_t\right)  $ of original equation (\ref{orginal1}) are given. Finally, we construct the stopping time and obtain appropriate control of $ u_{t}^{\epsilon}-\bar{u}_t $ before and after the stopping time respectively.

	\subsection{The averaged equation}
	To define the averaged equation, we first consider the equation (\ref{en12}) associated to the fast equation. 
	
	Under the assumptions (A1)-(A5), it is easy to prove that the equation (\ref{en12}) has a unique strong solution $ v_t^{s,x,y}, $ which is a time homogeneous Markov process. 
	Moreover, use the same argument as \cite{liu2020averaging}, there exists some constant $ \beta _{p}^{*}>0 $ such that the following estimates hold and we will not give a detailed proof here: 
	\begin{align}\label{en13}
	\mathbb{E}\left| v_{t}^{s,x,y} \right|^p\le C_{p,T} \left( 1+\left| x \right|^p \right) +e^{-\beta _{p}^{*} t}\left| y \right|^p,
	\end{align}
	and
	\begin{align}\label{en14}
	\mathbb{E}\left| v_{t}^{s,x,y_1}-v_{t}^{s,x,y_2} \right|^2\le e^{-\beta _1t}\left| y_1-y_2 \right|^2,
	\end{align}
	and
	\begin{eqnarray}\label{en17}
	\mathbb{E}\left| v_{t}^{s_1,x_1,y}-v_{t}^{s_2,x_2,y} \right|^2\le C_T \big( \left| s_1-s_2 \right|^{2\iota}+\left| x_1-x_2 \right|^2 \big)\big(1+\left| x_1\right|^{2\alpha_1} +\left| x_2\right|^{2\alpha_1\vee 2\alpha_2} +\left| y\right|^{2\alpha_2}  \big).
	\end{eqnarray} 
	
	Let $ \left\lbrace P_t^{s,x}\right\rbrace _{t\geq 0} $ be the transition semigroup of $ \left\lbrace v_t^{s,x,y}\right\rbrace_{t\geq 0},  $ that is 
	\begin{align}
	P_{t}^{s,x}\varphi \left( y \right) :=\mathbb{E}\varphi \left( v_{t}^{s,x,y} \right), \quad s>0,\ y\in \mathbb{R}^m,
	\end{align}
	where $ \varphi:\mathbb{R}^m\rightarrow \mathbb{R} $ is a bounded measurable function. 
	
	Then, we can establish the following crucial lemma:
	\begin{lem}\label{lem4.1}
		Assume that the conditions (A1)-(A5) hold. Then, for any fixed $ s>0  $ and $  x\in\mathbb{R}^n,    $  there
		exists a unique invariant measure  $ \mu ^{s,x}   $ for the equation (\ref{en12}), and 
		\begin{align}\label{en15}
		\int_{\mathbb{R}^m}{\left| z \right|^p\mu ^{s,x}\left( dz \right) \le C_{p,T} \left( 1+\left| x \right|^p \right)}.
		\end{align}
		Moreover, for any $ t>0  $ and $ y \in \mathbb{R}^m, $ we obtain
		\begin{align}\label{en16}
		\Big| \mathbb{E}b_1\left( s,x,v_{t}^{s,x,y} \right) -\int_{\mathbb{R}^m}{b_1\left( s,x,z \right) \mu ^{s,x}\left( dz \right)} \Big|\le C_Te^{-\frac{\beta _1}{2}t} \big( 1+\left| x \right|^{2\left( \theta _2\lor \theta _3\lor 1 \right)}+\left| y \right|^{2\left( \theta _3\lor 1 \right)}  \big).
		\end{align}
	\end{lem}
	\para{Proof:} The detailed proof will be given in the Appendix.

	Further, by the unique invariant measure $ \mu ^{s,x}, $ the averaged equation can be defined  as (\ref{en45}). Moreover, we can give some properties of the averaged coefficient $ \bar{b}_1, $ where the detailed proof of \lemref{lemm6.1} can be found in the Appendix.
	\begin{lem}\label{lemm6.1}
		Assume that the conditions {\rm (A1)-(A5)} hold. Then, for any $ t\geq 0 $ and $x \in \mathbb{R}^n,$ we have
		\begin{align}\label{en20}
		\left| \bar{b}_1\left( t,x \right) \right|\le C\left( 1+\left| x \right| \right).
		\end{align}
		Moreover, for any $ s_1, s_2\in [0,T], \ R\in\mathbb{R} $ and $ x_i\in \mathbb{R}^n   $ with  $ \left|x_i \right|\leq R, $ we have
		\begin{eqnarray}\label{en19}
		\left| \bar{b}_1\left( s_1,x_1 \right) -\bar{b}_1\left( s_2,x_2 \right) \right|\le C_{R,T}\big( \left| s_1-s_2 \right|^\iota+\left| s_1-s_2 \right|^\kappa+\left| x_1-x_2 \right| \big).
		%		\left( 1+\left| x_1 \right|^{\varrho_1}+\left| x_2 \right|^{\varrho_2} \right),
		\end{eqnarray}
		%	where
		%	$ \varrho_1:=\theta _1\lor \theta _2\lor 2\left( \alpha _1\lor \theta _3 \right) $ and $ \varrho_2:=2\left( \alpha _1\lor \alpha _2\lor \theta _2\lor \theta _3 \right). $
	\end{lem}
	
	Thanks to the assumptions and \lemref{lemm6.1}, by proceeding as \thmref{thm3.1} and \lemref{lem3.3}, it is easy to get that there exists a unique strong solution $ \bar{u}_t $ to equation (\ref{en45}) and we have
	\begin{align}\label{en46}
	\mathbb{E}\left\|  \bar{u} \right\|  _{\alpha,\infty  }^{p}\le C_{\alpha,p,x, T},\quad p\geq 1.  
	\end{align}

	\subsection{Some a-priori estimates}\label{sec-6} 
	To prove the \thmref{thm2.1}, some a-priori estimates for the solution $ \left(  u^\epsilon_t,  v^  \epsilon_t\right)  $ of original equation (\ref{orginal1}) need to be given at first.
	\begin{lem}\label{lem3.1}
		Assume that the conditions (A1)-(A5) hold. Then, for any $ \epsilon \in \left( 0,1\right], \ \alpha \in \left(  1-H, \frac{1}{2} \land\beta  \land  \frac{\gamma}{2} \right), \ p\geq 1 $ and $ t\in [0,T], $ we have
		\begin{align}\label{en22}
		\mathbb{E}\left\|  u^{\epsilon} \right\|  _{\alpha,\infty  }^{p}\le C_{\alpha,p,x, T}  \quad \text{and}\quad\mathbb{E}\left| v_{t}^{\epsilon} \right|^2\le C_{\alpha, x,y, T}.  
		\end{align}
		Moreover, for any $ h\in (0,1], $ it yields
		\begin{align}\label{en11}
		\mathbb{E}\left| u_{t+h}^{\epsilon}-u_{t}^{\epsilon} \right|^2\le C_{\alpha,x,y, T} h.
		\end{align}
	\end{lem}
	\para{Proof:}  
	According to the assumptions and use the same argument as \lemref{lem3.3} and \cite[Lemma 3.1]{liu2020averaging}, it is easy to get that equation (\ref{en22}) hold. Moreover, using the It\^{o} isometry for Brownian motion term and by 
	proceeding as Lemma 4.2 in \cite{pei2020averaging}, we also can establish the equation (\ref{en11}). Here, we omit the detailed proof. \qed

	Then, inspired by Khasminskii's idea in \cite{khas1968on},  
	for any $ \epsilon>0, $ the interval $ [0,T] $ is divided into subintervals of size $ \delta_\epsilon>0, $ where $ \delta_\epsilon $ is a fixed number depending on $ \epsilon. $ 
	Now, we construct a process $ \hat{v}  ^\epsilon_t $ with initial value $ \hat{v}  ^\epsilon_0= v^  \epsilon_0=y, $ and for $ t\in \left[k\delta_\epsilon, \text{min}\left\lbrace (k+1)\delta_\epsilon,T \right\rbrace  \right],  $  we have
	\begin{align}
	\hat{v}  _{t}^{\epsilon}=\hat{v}  _{k\delta _{\epsilon}}^{\epsilon}+\frac{1}{\epsilon}\int_{k\delta _{\epsilon}}^t{b_2\big( k\delta _{\epsilon},u_{k\delta _{\epsilon}}^{\epsilon},\hat{v}  _{r}^{\epsilon} \big) dr}+\frac{1}{\sqrt{\epsilon}}\int_{k\delta _{\epsilon}}^t{f_2\big( k\delta _{\epsilon},u_{k\delta _{\epsilon}}^{\epsilon},\hat{v}  _{r}^{\epsilon} \big) dW_{r}^{2}},
	\end{align}
	i.e.,
	\begin{align}
	\hat{v}  _{t}^{\epsilon}=y+\frac{1}{\epsilon}\int_0^t{b_2\big( r\left( \delta _{\epsilon} \right) ,u_{r\left( \delta _{\epsilon} \right)}^{\epsilon},\hat{v}  _{r}^{\epsilon} \big) dr}+\frac{1}{\sqrt{\epsilon}}\int_0^t{f_2\big( r\left( \delta _{\epsilon} \right) ,u_{r\left( \delta _{\epsilon} \right)}^{\epsilon},\hat{v}  _{r}^{\epsilon} \big) dW_{r}^{2}},
	\end{align}
	where $ r(\delta_{\epsilon})=\lfloor r/\delta_{\epsilon}\rfloor\delta_{\epsilon} $ is the nearest breakpoint preceding $ r. $ 
	By the construction of $ \hat{v}  ^\epsilon_t, $ we have an estimate analogous to \lemref{lem3.1} hold, i.e., for any $ t\in [0,T], $ we have
	\begin{align}\label{en24}
	\mathbb{E}\big| \hat{v}  _{t}^{\epsilon} \big|^2\le C_{\alpha ,x,y, T}.
	\end{align}
	Moreover, thanks to the assumptions and \lemref{lem3.1},  by proceeding as Lemma 3.4 in \cite{liu2020averaging}, it is easy to get that
	\begin{align}\label{en23}
	\mathbb{E}\big| v_{t}^{\epsilon}-\hat{v}  _{t}^{\epsilon} \big|^2\le C_{\alpha ,x,y, T}\delta _{\epsilon}^{2\iota \land 1}. 
	\end{align}
	
	\subsection{The proof of  \thmref{thm2.1}}
	Now, we construct the following stopping time $ \tau _{R}^{\epsilon}  $ for each $ R\in \mathbb{R}: $
	\begin{align}
	\tau _{R}^{\epsilon}:=\text{inf}\left\{ t\geq 0:\left\|  u_{t}^{\epsilon} \right\|  _{\alpha}+\left\|  \bar{u}_t \right\|  _{\alpha}+\varLambda _{\alpha}^{0,t}\left( B^H \right) \geq R \right\} 
	.\nonumber
	\end{align}
	Moreover, due to the trajectories of  $ u^\epsilon_t $ and $ \bar{u}_t $ are $ \eta $-H\"{o}lder continuous for all $ \eta< 1/2. $ 
	As the proof in \thmref{thm3.1}, let $ \eta\in ( \frac{\alpha}{\gamma}, \frac{1}{2})   $ and consider the following set $ \Omega_M\subset\Omega $ with $ M\in\mathbb{N},  $ such that
	\begin{align}
	\Omega_M:=\big\lbrace \omega \in \Omega:\left\| u^\epsilon\right\| _\eta\leq M\  \text{and} \ \left\| \bar{u}\right\| _\eta\leq M \big\rbrace. \nonumber
	\end{align}
	It is clear that $ \Omega_M\nearrow\Omega. $ 
	
	First, we estimate the error of $ u_{t}^{\epsilon}-\bar{u}_t $ before a
	stopping time: 
	\begin{lem} 
		Assume that the conditions {\rm (A1)-(A5)} hold. Then, for any $ \alpha \in \left(  1-H, \frac{1}{2} \land\beta\right.   $  $ \left. \land  \frac{\gamma}{2} \right),  $ we have
		\begin{align}\label{en49}
		\underset{t\in \left[ 0,T \right]}{\sup}\mathbb{E}\left[ \left\|  u_{t}^{\epsilon}-\bar{u}_t \right\|  _{\alpha}^{2}\mathbf{1}_{\varOmega _M\cap \{T\leq \tau _{R}^{\epsilon}\}} \right] \leq C_{\alpha ,x,y, R,M,T}\big(\delta _{\epsilon}^{2\kappa \land 2\iota \land \left( {1}/{2} -\alpha\right)}+ {\epsilon}/{\delta _{\epsilon}}\big). 
		\end{align}
	\end{lem}
	\para{Proof:}  
	%Note that
	%\begin{align}
	%u_{t}^{\epsilon}-\bar{u}_t&=\int_0^t{\left( b_1\left( r,u_{r}^{\epsilon},v_{r}^{\epsilon} \right) -\bar{b}_1\left( r,\bar{u}_r \right) \right) dr}+\int_0^t{\left( f_1\left( r,u_{r}^{\epsilon} \right) -f_1\left( r,\bar{u}_r \right) \right) dW_{r}^{1}}\cr
	%&\quad+\int_0^t{\left( g_1\left( r,u_{r}^{\epsilon} \right) -g_1\left( r,\bar{u}_r \right) \right) dB_{r}^{H}}.\nonumber
	%\end{align}
	It is easy to know that
	\begin{align}\label{en34}
	\mathbb{E}\left[ \left\|  u_{t}^{\epsilon}-\bar{u}_t \right\|  _{\alpha}^{2}\mathbf{1}_{\varOmega _M\cap \{T\leq \tau _{R}^{\epsilon}\}} \right] &\le 3\mathbb{E}\Big[ \Big\| \int_0^t{\left( b_1\left( r,u_{r}^{\epsilon},v_{r}^{\epsilon} \right) -\bar{b}_1\left( r,\bar{u}_r \right) \right) dr} \Big\| _{\alpha}^{2}\mathbf{1}_{\varOmega _M\cap \{T\leq \tau _{R}^{\epsilon}\}} \Big] \cr
	&\quad+3\mathbb{E}\Big[ \Big\| \int_0^t{\left( f_1\left( r,u_{r}^{\epsilon} \right) -f_1\left( r,\bar{u}_r \right) \right) dW_{r}^{1}} \Big\| _{\alpha}^{2}\mathbf{1}_{\varOmega _M\cap \{T\leq \tau _{R}^{\epsilon}\}} \Big] \cr
	&\quad+3\mathbb{E}\Big[ \Big\| \int_0^t{\left( g_1\left( r,u_{r}^{\epsilon} \right) -g_1\left( r,\bar{u}_r \right) \right) dB_{r}^{H}} \Big\| _{\alpha}^{2}\mathbf{1}_{\varOmega _M\cap \{T\leq \tau _{R}^{\epsilon}\}} \Big] \cr
	&:=3 \mathcal{I}_{t}^{1}+3\mathcal{I}_{t}^{2}+3\mathcal{I}_{t}^{3}.
	\end{align}
	
	For $ \mathcal{I}_t^2  $ and $ \mathcal{I}_t^3,  $ if $ \omega\in \Omega_n $ and $ \eta\in ( \frac{\alpha}{\gamma}, \frac{1}{2}), $ by proceeding as (\ref{en525}), use the conclusions of Proposition 3.6 and Proposition 3.9 in \cite{Guerra2008Stochastic}, it is easy to get that
	\begin{align}\label{en32}
	\mathcal{I}_{t}^{2}+\mathcal{I}_{t}^{3}\le  
	C_{\alpha ,R,M,T}  \int_0^t{\left( t-r \right) ^{-\frac{1}{2}-\alpha} r^{-\frac{1}{2}-\alpha} \mathbb{E}\left[ \left\|  u_{r}^{\epsilon}-\bar{u}_r \right\|  _{\alpha}^{2}\mathbf{1}_{\varOmega _M\cap \{T\leq \tau _{R}^{\epsilon}\}} \right] dr}.
	\end{align}

	For $ \mathcal{I}_t^1, $ we have
	\begin{align}
	\mathcal{I}_{t}^{1}
	%&=\mathbb{E}\Big[ \Big\|  \int_0^t{b_1\left( r,u_{r}^{\epsilon},v_{r}^{\epsilon} \right) -b_1\left( r,\bar{u}_r \right) dr} \Big\|  _{\alpha}^{2}\mathbf{1}_{\varOmega _M} \Big]\cr 
	&\le 4\mathbb{E}\Big[ \Big\|  \int_0^t{b_1\left( r,u_{r}^{\epsilon},v_{r}^{\epsilon} \right) -b_1\big( r\left( \delta _{\epsilon} \right) ,u_{r\left( \delta _{\epsilon} \right)}^{\epsilon},\hat{v}  _{r}^{\epsilon} \big) dr} \Big\|  _{\alpha}^{2}\mathbf{1}_{\varOmega _M\cap \{T\leq \tau _{R}^{\epsilon}\}} \Big]\cr
	&\quad+4\mathbb{E}\Big[\Big\|  \int_0^t{b_1\big( r\left( \delta _{\epsilon} \right) ,u_{r\left( \delta _{\epsilon} \right)}^{\epsilon},\hat{v}  _{r}^{\epsilon} \big) -\bar{b}_1\left( r\left( \delta _{\epsilon} \right) ,u_{r\left( \delta _{\epsilon} \right)}^{\epsilon} \right) dr} \Big\|  _{\alpha}^{2}\mathbf{1}_{\varOmega _M\cap \{T\leq \tau _{R}^{\epsilon}\}} \Big]\cr
	&\quad+4\mathbb{E}\Big[ \Big\|  \int_0^t{\bar{b}_1\left( r\left( \delta _{\epsilon} \right) ,u_{r\left( \delta _{\epsilon} \right)}^{\epsilon} \right) -\bar{b}_1\left( r,u_{r}^{\epsilon} \right) dr} \Big\|  _{\alpha}^{2}\mathbf{1}_{\varOmega _M\cap \{T\leq \tau _{R}^{\epsilon}\}}  \Big]\cr
	&\quad+4\mathbb{E}\Big[ \Big\|  \int_0^t{\bar{b}_1\left( r,u_{r}^{\epsilon} \right) -\bar{b}_1\left( r,\bar{u}_r \right) dr} \Big\|  _{\alpha}^{2}\mathbf{1}_{\varOmega _M\cap \{T\leq \tau _{R}^{\epsilon}\}} \Big]\cr 
	&:=4 \mathcal{J}_{t}^{1}+4\mathcal{J}_{t}^{2}+4\mathcal{J}_{t}^{3}+4\mathcal{J}_{t}^{4} .
	\end{align}
	Thanks to \lemref{lem3.1} and equation (\ref{en23}), we can get
	\begin{align}
	&\qquad\mathbb{E}\Big[ \Big| \int_s^t{b_1\left( r,u_{r}^{\epsilon},v_{r}^{\epsilon} \right) -b_1\big( r\left( \delta _{\epsilon} \right) ,u_{r\left( \delta _{\epsilon} \right)}^{\epsilon},\hat{v}  _{r}^{\epsilon} \big) dr} \Big|^2\mathbf{1}_{\varOmega _M\cap \{T\leq \tau _{R}^{\epsilon}\}} \Big] \cr
	&\le \mathbb{E}\Big[ \Big| \int_s^t{b_1\left( r,u_{r}^{\epsilon},v_{r}^{\epsilon} \right) -b_1\left( r\left( \delta _{\epsilon} \right) ,u_{r}^{\epsilon},v_{r}^{\epsilon} \right) dr} \Big|^2\mathbf{1}_{\varOmega _M\cap \{T\leq \tau _{R}^{\epsilon}\}} \Big] \cr
	&\quad+\mathbb{E}\Big[ \Big| \int_s^t{b_1\left( r\left( \delta _{\epsilon} \right) ,u_{r}^{\epsilon},v_{r}^{\epsilon} \right) -b_1\left( r\left( \delta _{\epsilon} \right) ,u_{r\left( \delta _{\epsilon} \right)}^{\epsilon},v_{r}^{\epsilon} \right) dr} \Big|^2\mathbf{1}_{\varOmega _M\cap \{T\leq \tau _{R}^{\epsilon}\}} \Big] \cr
	&\quad+\mathbb{E}\Big[ \Big| \int_s^t{b_1\left( r\left( \delta _{\epsilon} \right) ,u_{r\left( \delta _{\epsilon} \right)}^{\epsilon},v_{r}^{\epsilon} \right) -b_1\big( r\left( \delta _{\epsilon} \right) ,u_{r\left( \delta _{\epsilon} \right)}^{\epsilon},\hat{v}  _{r}^{\epsilon} \big) dr} \Big|^2\mathbf{1}_{\varOmega _M\cap \{T\leq \tau _{R}^{\epsilon}\}} \Big] \cr
	&\le C_T\int_s^t{\left| r-r\left( \delta _{\epsilon} \right) \right|^{2\kappa}\mathbb{E}\big[ \big( 1+\left| u_{r}^{\epsilon} \right|^{2\theta _2}+\left| v_{r}^{\epsilon} \right|^{2\theta _3} \big) \mathbf{1}_{\varOmega _M\cap \{T\leq \tau _{R}^{\epsilon}\}} \big]}dr\cr
	&\quad+C_{R,T}\int_s^t{\mathbb{E}\big[ \left| u_{r}^{\epsilon}-u_{r\left( \delta _{\epsilon} \right)}^{\epsilon} \right|^2\mathbf{1}_{\varOmega _M\cap \{T\leq \tau _{R}^{\epsilon}\}} \big]}dr\int_s^t{\mathbb{E}\big[ \big( 1+\left| v_{r}^{\epsilon} \right|^{2\theta _1} \big) \mathbf{1}_{\varOmega _M\cap \{T\leq \tau _{R}^{\epsilon}\}} \big]}dr\cr
	&\quad+C_T\int_s^t{\mathbb{E}\big[ \big| v_{r}^{\epsilon}-\hat{v}  _{r}^{\epsilon} \big|^2\mathbf{1}_{\varOmega _M\cap \{T\leq \tau _{R}^{\epsilon}\}} \big]}dr\int_s^t{\mathbb{E}\big[ \big( 1+\left| u_{r}^{\epsilon} \right|^{2\theta _2}+\left| v_{r}^{\epsilon} \right|^{2\theta _3}+\big| \hat{v}  _{r}^{\epsilon} \big|^{2\theta _3} \big) \mathbf{1}_{\varOmega _M\cap \{T\leq \tau _{R}^{\epsilon}\}} \big]}dr\cr
	&\le C_{\alpha ,x,y,R,T}\left( \delta _{\epsilon}^{2\kappa}+\delta _{\epsilon}+\delta _{\epsilon}^{2\iota \land 1} \right)\left( t-s\right).
	\end{align}
	Hence 
	\begin{align}\label{en29}
	\mathcal{J}_{t}^{1}
	%&=\mathbb{E}\Big[ \Big| \int_0^t{b_1\left( r,u_{r}^{\epsilon},v_{r}^{\epsilon} \right) -b_1\big( r\left( \delta _{\epsilon} \right) ,u_{r\left( \delta _{\epsilon} \right)}^{\epsilon},\hat{v}  _{r}^{\epsilon} \big) dr} \Big|\cr
	%&\qquad\qquad\qquad\qquad +\int_0^t{\frac{\left| \int_r^t{b_1\left( s,u_{s}^{\epsilon},v_{s}^{\epsilon} \right) -b_1\big( s\left( \delta _{\epsilon} \right) ,u_{s\left( \delta _{\epsilon} \right)}^{\epsilon},\hat{v}  _{s}^{\epsilon} \big) ds} \right|}{\left( t-r \right) ^{\alpha +1}}dr} \Big] ^2\cr
	&\le 2\mathbb{E}\Big[ \Big| \int_0^t{b_1\left( r,u_{r}^{\epsilon},v_{r}^{\epsilon} \right) -b_1\big( r\left( \delta _{\epsilon} \right) ,u_{r\left( \delta _{\epsilon} \right)}^{\epsilon},\hat{v}  _{r}^{\epsilon} \big) dr} \Big|^2\mathbf{1}_{\varOmega _M\cap \{T\leq \tau _{R}^{\epsilon}\}}  \Big]\cr
	&\quad+2\mathbb{E}\Big[\Big| \int_0^t{\left( t-r \right) ^{-\alpha -1}\Big| \int_r^t{b_1\left( s,u_{s}^{\epsilon},v_{s}^{\epsilon} \right) -b_1\big( s\left( \delta _{\epsilon} \right) ,u_{s\left( \delta _{\epsilon} \right)}^{\epsilon},\hat{v}  _{s}^{\epsilon} \big) ds} \Big|dr} \Big|^2\mathbf{1}_{\varOmega _M\cap \{T\leq \tau _{R}^{\epsilon}\}}  \Big]\cr
	&\le C_{\alpha ,x,y,R,T}\left( \delta _{\epsilon}^{2\kappa}+\delta _{\epsilon}+\delta _{\epsilon}^{2\iota \land 1} \right)+ C\int_0^t{\left( t-r \right) ^{-\alpha -\frac{1}{2}}dr}\cr 
	&\qquad\times\int_0^t{\left( t-r \right) ^{-\alpha -\frac{3}{2}}\mathbb{E}\Big[ \Big| \int_r^t{b_1\left( s,u_{s}^{\epsilon},v_{s}^{\epsilon} \right) -b_1\big( s\left( \delta _{\epsilon} \right) ,u_{s\left( \delta _{\epsilon} \right)}^{\epsilon},\hat{v}  _{s}^{\epsilon} \big) ds} \Big|^2\mathbf{1}_{\varOmega _M\cap \{T\leq \tau _{R}^{\epsilon}\}} \Big]}dr
	\cr
	&\le C_{\alpha ,x,y,R,T}\left( \delta _{\epsilon}^{2\kappa}+\delta _{\epsilon}+\delta _{\epsilon}^{2\iota \land 1} \right) \le C_{\alpha ,x,y,R,T}\delta _{\epsilon}^{2\kappa \land 2\iota \land 1}.
	\end{align}
	Then, due to (\ref{en19}) and Proposition 3.4 in \cite{Guerra2008Stochastic}, by proceeding as (\ref{en29}), it is easy to get that
	\begin{eqnarray}\label{en30}
	\mathcal{J}_{t}^{3}+\mathcal{J}_{t}^{4}\le C_{\alpha ,R,x,y, T}\delta _{\epsilon}^{2\kappa \land 2\iota \land 1}+ C_{\alpha ,R,T}\int_0^t{\left( t-r \right) ^{-2\alpha}\mathbb{E}\left[ \left\|  u_{r}^{\epsilon}-\bar{u}_r \right\|  _{\alpha}^{2}\mathbf{1}_{\varOmega _M\cap \{T\leq \tau _{R}^{\epsilon}\}} \right] dr}.
	\end{eqnarray}
	
	In order to prove the approximation result of the equation (\ref{en49}), we must estimate $ \mathcal{J}_{t}^{2}. $ Therefore, we establish the following crucial lemma:
	\begin{lem}\label{lem5.2}
		Assume that the conditions {\rm (A1)-(A5)} hold. Then, for any $ s,\varLambda>0, x  \in \mathbb{R}^n $ and $ y\in\mathbb{R}^m, $ we have
		\begin{align}\label{en21}
		\mathbb{E}\Big| \frac{1}{\varLambda}\int_0^{\varLambda}{b_1\left( s,x,v_{r}^{s,x,y} \right) dr}-\bar{b}_1\left( s,x \right)   \Big|^2\le \frac{C}{\varLambda}\left( 1+\left| x \right|^{4\left( \theta _2\lor \theta _3\lor 1 \right)}+\left| y \right|^{4\left( \theta _3\lor 1 \right)} \right). 
		\end{align}
	\end{lem}
	\para{Proof:}
	Using the Markov property of $ \left\lbrace  v^  {s,x,y}_t\right\rbrace_{t\geq 0},  $ we have
	\begin{align}
	&\qquad\mathbb{E}\Big| \frac{1}{\varLambda}\int_0^\varLambda{b_1\left( s,x,v_{r}^{s,x,y} \right) dr}-\bar{b}_1\left( s,x \right) \Big|^2\cr
	&=\frac{2}{\varLambda^2}\int_0^\varLambda{\int_r^\varLambda{\mathbb{E}\big[  \left( b_1\left( s,x,v_{r}^{s,x,y} \right) -\bar{b}_1\left( s,x \right) \right) \left( b_1\left( s,x,v_{\sigma}^{s,x,y} \right)  -\bar{b}_1\left( s,x \right) \right)  \big] d\sigma dr}}\cr
	&=\frac{2}{\varLambda^2}\int_0^\varLambda{\int_r^\varLambda{\mathbb{E}\big[ \left( b_1\left( s,x,v_{r}^{s,x,y} \right) -\bar{b}_1\left( s,x \right) \right) P_{\sigma -r}^{s,x}\left( b_1\left( s,x,v_{r}^{s,x,y} \right)  -\bar{b}_1\left( s,x \right) \right) \big] d\sigma dr}}\cr
	&\le \frac{2}{\varLambda^2}\int_0^\varLambda{\int_r^\varLambda{\big[ \mathbb{E}\left| b_1\left( s,x,v_{r}^{s,x,y} \right) -\bar{b}_1\left( s,x \right) \right|^2   \mathbb{E}\left| P_{\sigma -r}^{s,x}\left( b_1\left( s,x,v_{r}^{s,x,y} \right) -\bar{b}_1\left( s,x \right) \right) \right|^2 \big] ^{\frac{1}{2}}d\sigma dr}}.\nonumber
	\end{align}
	By assumption (A1) and thanks to 
	(\ref{en13}), (\ref{en20}), we can get
	\begin{align}
	\mathbb{E}\left| b_1\left( s,x,v_{r}^{s,x,y} \right) -\overline{b}_1\left( s,x \right) \right|^2&\le 2\mathbb{E}\left| b_1\left( s,x,v_{r}^{s,x,y} \right) \right|^2+2\mathbb{E}\left| \overline{b}_1\left( s,x \right) \right|^2\cr
	&\le C \left( 1+\left| x \right|^2+\mathbb{E}\left| v_{r}^{s,x,y} \right|^2 \right) \cr
	&\le C \left( 1+\left| x \right|^2 \right) +e^{-\beta _2^* t}\left| y \right|^2.\nonumber
	\end{align}
	Then, according to the equation (\ref{en16}) and (\ref{en13}), we obtain
	\begin{align}
	\mathbb{E}\left| P_{\sigma -r}^{s,x}\left( b_1\left( s,x,v_{r}^{s,x,y} \right) dr-\bar{b}_1\left( s,x \right) \right) \right|^2&=\mathbb{E}\Big| \mathbb{E}b_1\big( s,x,v_{\sigma -r}^{s,x,v_{r}^{s,x,y}} \big) -\int_{\mathbb{R}^m}{b_1\left( s,x,z \right) \mu ^{s,x}\left( dz \right)} \Big|^2\cr
	&\le C_Te^{-\frac{\beta _1}{2}\left( \sigma -r \right)}\mathbb{E}\big| 1+\left| x \right|^{2\left( \theta _2\lor \theta _3\lor 1 \right)}+\mathbb{E}\left| v_{r}^{s,x,y} \right|^{2\left( \theta _3\lor 1 \right)} \big|^2\cr
	&\le C_Te^{-\frac{\beta _1}{2}\left( \sigma -r \right)}\big( 1+\left| x \right|^{4\left( \theta _2\lor \theta _3\lor 1 \right)}+\left| y \right|^{4\left( \theta _3\lor 1 \right)} \big).\nonumber
	\end{align}
	Hence
	\begin{align}
	\mathbb{E}\Big| \frac{1}{\varLambda}\int_0^{\varLambda}{b_1\left( s,x,v_{r}^{s,x,y} \right) dr}-\bar{b}_1\left( s,x \right) \Big|^2&\le \frac{C}{\varLambda ^2}\big( 1+\left| x \right|^{4\left( \theta _2\lor \theta _3\lor 1 \right)}+\left| y \right|^{4\left( \theta _3\lor 1 \right)} \big)\int_0^{\varLambda}{\int_r^{\varLambda}{e^{-\frac{\beta _1}{2}\left( \sigma -r \right)}d\sigma dr}}\cr
	&\le \frac{C}{\varLambda}\big( 1+\left| x \right|^{4\left( \theta _2\lor \theta _3\lor 1 \right)}+\left| y \right|^{4\left( \theta _3\lor 1 \right)} \big).\nonumber
	\end{align}
	The proof is completed.\qed
	\begin{lem}\label{lem6.2}
		Assume that the conditions {\rm (A1)-(A5)} hold. Then, for any $ t\in [0,T], $ we have
		\begin{align}\label{en28}
		\mathbb{E}\Big\|  \int_0^t{b_1\big( r\left( \delta _{\epsilon} \right) ,u_{r\left( \delta _{\epsilon} \right)}^{\epsilon},\hat{v}  _{r}^{\epsilon} \big) -\bar{b}_1\left( r\left( \delta _{\epsilon} \right) ,u_{r\left( \delta _{\epsilon} \right)}^{\epsilon} \right) dr} \Big\| _{\alpha}^{2}\leq C_{\alpha ,x,y, T}\big( \frac{\epsilon}{\delta _{\epsilon}}+  \delta _{\epsilon}^{\frac{1}{2}-\alpha}\big).
		\end{align}
	\end{lem}
	\para{Proof:} By elementary inequality, we have
	\begin{align}\label{en25}
	&\qquad\mathbb{E}\Big\|  \int_0^t{b_1\big( r\left( \delta _{\epsilon} \right) ,u_{r\left( \delta _{\epsilon} \right)}^{\epsilon},\hat{v}  _{r}^{\epsilon} \big) -\bar{b}_1\left( r\left( \delta _{\epsilon} \right) ,u_{r\left( \delta _{\epsilon} \right)}^{\epsilon} \right) dr} \Big\| _{\alpha}^{2}\cr
	&\le 3\mathbb{E}\Big| \sum_{k=0}^{\lfloor t/\delta _{\epsilon} \rfloor -1}{\int_{k\delta _{\epsilon}}^{\left( k+1 \right) \delta _{\epsilon}}{b_1\big( k\delta _{\epsilon},u_{k\delta _{\epsilon}}^{\epsilon},\hat{v}  _{r}^{\epsilon} \big) -\bar{b}_1\left( k\delta _{\epsilon},u_{k\delta _{\epsilon}}^{\epsilon} \right) dr}} \Big|^2\cr
	&\quad+3\mathbb{E}\Big| \int_{\lfloor t/\delta _{\epsilon} \rfloor \delta _{\epsilon}}^t{b_1\big( r\left( \delta _{\epsilon} \right) ,u_{r\left( \delta _{\epsilon} \right)}^{\epsilon},\hat{v}  _{r}^{\epsilon} \big) -\bar{b}_1\left( r\left( \delta _{\epsilon} \right) ,u_{r\left( \delta _{\epsilon} \right)}^{\epsilon} \right) dr} \Big|^2\cr
	&\quad+3\mathbb{E}\Big| \int_0^t{\left( t-r \right) ^{-\alpha -1}\Big| \int_r^t{b_1\big( s\left( \delta _{\epsilon} \right) ,u_{s\left( \delta _{\epsilon} \right)}^{\epsilon},\hat{v}  _{s}^{\epsilon} \big) -\bar{b}_1\left( s\left( \delta _{\epsilon} \right) ,u_{s\left( \delta _{\epsilon} \right)}^{\epsilon} \right) ds} \Big|dr} \Big|^2\cr
	&:=
	\mathcal{K}_{t}^{1}+\mathcal{K}_{t}^{2}+\mathcal{K}_{t}^{3}.
	\end{align}
	According to \lemref{lem5.2} and \lemref{lem3.1},  we can get 
	\begin{align}\label{en26}
	\mathcal{K}_{t}^{1}&\le 3\lfloor t/\delta _{\epsilon} \rfloor \sum_{k=0}^{\lfloor t/\delta _{\epsilon} \rfloor -1}{\mathbb{E}\Big| \int_{k\delta _{\epsilon}}^{\left( k+1 \right) \delta _{\epsilon}}{b_1\big( k\delta _{\epsilon},u_{k\delta _{\epsilon}}^{\epsilon},\hat{v}  _{r}^{\epsilon} \big) -\bar{b}_1\left( k\delta _{\epsilon},u_{k\delta _{\epsilon}}^{\epsilon} \right) dr} \Big|}^2\cr
	%&\le \frac{C_T}{\delta _{\epsilon}^{2}}\underset{0\le k\le \lfloor t/\delta _{\epsilon} \rfloor -1}{\max}\mathbb{E}\Big| \int_{k\delta _{\epsilon}}^{\left( k+1 \right) \delta _{\epsilon}}{b_1\big( k\delta _{\epsilon},u_{k\delta _{\epsilon}}^{\epsilon},\hat{v}  _{r}^{\epsilon} \big) -\bar{b}_1\left( k\delta _{\epsilon},u_{k\delta _{\epsilon}}^{\epsilon} \right) dr} \Big|^2\cr
	&\le \frac{C_T}{\delta _{\epsilon}^{2}}\underset{0\le k\le \lfloor t/\delta _{\epsilon} \rfloor -1}{\max}\mathbb{E}\Big| \int_0^{\delta _{\epsilon}}{b_1\big( k\delta _{\epsilon},u_{k\delta _{\epsilon}}^{\epsilon},\hat{v}  _{k\delta _{\epsilon}+r}^{\epsilon} \big) -\bar{b}_1\left( k\delta _{\epsilon},u_{k\delta _{\epsilon}}^{\epsilon} \right) dr} \Big|^2\cr
	&\le C_T\frac{\epsilon ^2}{\delta _{\epsilon}^{2}}\underset{0\le k\le \lfloor t/\delta _{\epsilon} \rfloor -1}{\max}\mathbb{E}\Big| \int_0^{ {\delta _{\epsilon}}/{\epsilon}}{b_1\big( k\delta _{\epsilon},u_{k\delta _{\epsilon}}^{\epsilon},\hat{v}  _{k\delta _{\epsilon}+\epsilon r}^{\epsilon} \big) -\bar{b}_1\left( k\delta _{\epsilon},u_{k\delta _{\epsilon}}^{\epsilon} \right) dr} \Big|^2\cr
	&\le C_T\underset{0\le k\le \lfloor t/\delta _{\epsilon} \rfloor -1}{\max}\mathbb{E}\Big| \frac{1}{ {\delta _{\epsilon}}/{\epsilon}}\int_0^{ {\delta _{\epsilon}}/{\epsilon}}{b_1\big( k\delta _{\epsilon},u_{k\delta _{\epsilon}}^{\epsilon}, {Y}_{r}^{k\delta _{\epsilon},u_{k\delta _{\epsilon}}^{\epsilon},v_{k\delta _{\epsilon}}^{\epsilon}} \big) -\bar{b}_1\left( k\delta _{\epsilon},u_{k\delta _{\epsilon}}^{\epsilon} \right) dr} \Big|^2\cr
	&\le C_{\alpha ,x,y, T}\frac{\epsilon}{\delta _{\epsilon}}. 
	\end{align}
	The last equation is thanks to the distribution of $ \big(u_{k\delta _{\epsilon}}^{\epsilon},\hat{v}  _{k\delta _{\epsilon}+  r}^{\epsilon} \big)   $ coincides with the distribution of $ \big( u_{k\delta _{\epsilon}}^{\epsilon}, {Y}_{r/\epsilon}^{k\delta _{\epsilon},u_{k\delta _{\epsilon}}^{\epsilon},v_{k\delta _{\epsilon}}^{\epsilon}} \big) $ in the interval $ r\in \left[0,\delta_{\epsilon} \right)   $ \cite{liu2020averaging}, where $ {Y}_{r}^{k\delta _{\epsilon},u_{k\delta _{\epsilon}}^{\epsilon},v_{k\delta _{\epsilon}}^{\epsilon}} $
	is the solution of the fast equation (\ref{en12}) by fixed $ k\delta _{\epsilon}> 0, $   frozen
	slow component $ u_{k\delta _{\epsilon}}^{\epsilon} $ and with initial datum $ v_{k\delta _{\epsilon}}^{\epsilon}, $ 
	and noise $ W^2 $ independent of
	both of them. 
	
	For $ \mathcal{K}_{t}^{2}, $ thanks to (\ref{en20}), (\ref{en22}) and (\ref{en24}),  we obtain 
	\begin{align}\label{en37}
	\mathcal{K}_{t}^{2}
	%&=3\mathbb{E}\Big| \int_{\lfloor t/\delta _{\epsilon} \rfloor \delta _{\epsilon}}^t{b_1\big( r\left( \delta _{\epsilon} \right) ,u_{r\left( \delta _{\epsilon} \right)}^{\epsilon},\hat{v}  _{r}^{\epsilon} \big) -\bar{b}_1\left( r\left( \delta _{\epsilon} \right) ,u_{r\left( \delta _{\epsilon} \right)}^{\epsilon} \right) dr} \Big|^2\cr
	&\le C\delta _{\epsilon}\int_{\lfloor t/\delta _{\epsilon} \rfloor \delta _{\epsilon}}^t{\mathbb{E}\big| b_1\big( r\left( \delta _{\epsilon} \right) ,u_{r\left( \delta _{\epsilon} \right)}^{\epsilon},\hat{v}  _{r}^{\epsilon} \big) \big|^2+\mathbb{E}\left| \bar{b}_1\left( r\left( \delta _{\epsilon} \right) ,u_{r\left( \delta _{\epsilon} \right)}^{\epsilon} \right) \right|^2dr} 
	%&\le C_T\delta _{\epsilon}\int_{\lfloor t/\delta _{\epsilon} \rfloor \delta _{\epsilon}}^t{\big( 1+\mathbb{E}\left| u_{r\left( \delta _{\epsilon} \right)}^{\epsilon} \right|^2+\mathbb{E}\big| \hat{v}  _{r}^{\epsilon} \big|^2 \big) dr}\cr
	&\le C_{\alpha ,x,y, T}\delta _{\epsilon}^{2}.
	\end{align}
	Then, using the H\"{o}lder inequality for $ \mathcal{K}_{t}^{3}, $ we have
	\begin{align}
	\mathcal{K}_{t}^{3}
	%&=3\mathbb{E}\Big| \int_0^t{\left( t-r \right) ^{-\alpha -1}\Big| \int_r^t{b_1\big( s\left( \delta _{\epsilon} \right) ,u_{s\left( \delta _{\epsilon} \right)}^{\epsilon},\hat{v}  _{s}^{\epsilon} \big) -\bar{b}_1\left( s\left( \delta _{\epsilon} \right) ,u_{s\left( \delta _{\epsilon} \right)}^{\epsilon} \right) ds} \Big|dr} \Big|^2\cr
	&\le 3\int_0^t{\left( t-r \right) ^{-\alpha -\frac{1}{2}}dr}\int_0^t{\left( t-r \right) ^{-\alpha -\frac{3}{2}}}\mathbb{E}\Big| \int_r^t{b_1\big( s\left( \delta _{\epsilon} \right) ,u_{s\left( \delta _{\epsilon} \right)}^{\epsilon},\hat{v}  _{s}^{\epsilon} \big) -\bar{b}_1\left( s\left( \delta _{\epsilon} \right) ,u_{s\left( \delta _{\epsilon} \right)}^{\epsilon} \right) ds} \Big|^2dr\cr
	&\le C_{\alpha ,T}\int_0^{\left( \lfloor \frac{t}{\delta _{\epsilon}} \rfloor -3 \right) \delta _{\epsilon}}{\left( t-r \right) ^{-\alpha -\frac{3}{2}}\mathbb{E}\Big| \int_r^t{b_1\big( s\left( \delta _{\epsilon} \right) ,u_{s\left( \delta _{\epsilon} \right)}^{\epsilon},\hat{v}  _{s}^{\epsilon} \big) -\bar{b}_1\left( s\left( \delta _{\epsilon} \right) ,u_{s\left( \delta _{\epsilon} \right)}^{\epsilon} \right) ds} \Big|^2dr}\cr
	&\quad+C_{\alpha ,T}\int_{\left( \lfloor \frac{t}{\delta _{\epsilon}} \rfloor -3 \right) \delta _{\epsilon}}^{t}{\left( t-r \right) ^{-\alpha -\frac{3}{2}}\mathbb{E}\Big| \int_r^t{b_1\big( s\left( \delta _{\epsilon} \right) ,u_{s\left( \delta _{\epsilon} \right)}^{\epsilon},\hat{v}  _{s}^{\epsilon} \big) -\bar{b}_1\left( s\left( \delta _{\epsilon} \right) ,u_{s\left( \delta _{\epsilon} \right)}^{\epsilon} \right) ds} \Big|^2dr}\cr
	%&\quad+C_{\alpha ,T}\int_{\lfloor \frac{t}{\delta _{\epsilon}} \rfloor \delta _{\epsilon}}^t{\left( t-r \right) ^{-\alpha -\frac{3}{2}}\mathbb{E}\Big| \int_r^t{b_1\left( s\left( \delta _{\epsilon} \right) ,u_{s\left( \delta _{\epsilon} \right)}^{\epsilon},v_{s}^{\epsilon} \right) -\bar{b}_1\left( s\left( \delta _{\epsilon} \right) ,u_{s\left( \delta _{\epsilon} \right)}^{\epsilon} \right) ds} \Big|^2dr}\cr
	&:=\mathcal{K}_{t}^{31}+\mathcal{K}_{t}^{32}.\nonumber
	\end{align}
	Using the same argument as the proof of (\ref{en26}), (\ref{en37})
	and the fact $ \lfloor\lambda_1\rfloor-\lfloor\lambda_2\rfloor\leq \lambda_1-\lambda_2+1.  $ Then, thanks to \lemref{lem5.2}, we can get 
	\begin{align}
	\mathcal{K}_{t}^{31}
	%&=C_{\alpha ,T}\int_0^{\left( \lfloor \frac{t}{\delta _{\epsilon}} \rfloor -3 \right) \delta _{\epsilon}}{\left( t-r \right) ^{-\alpha -\frac{3}{2}}\mathbb{E}\Big| \int_r^t{b_1\left( s\left( \delta _{\epsilon} \right) ,u_{s\left( \delta _{\epsilon} \right)}^{\epsilon},v_{s}^{\epsilon} \right) -\bar{b}_1\left( s\left( \delta _{\epsilon} \right) ,u_{s\left( \delta _{\epsilon} \right)}^{\epsilon} \right) ds} \Big|^2dr}\cr
	&\le C_{\alpha ,T}\int_0^{\left( \lfloor \frac{t}{\delta _{\epsilon}} \rfloor -3 \right) \delta _{\epsilon}}{\left( t-r \right) ^{-\alpha -\frac{3}{2}}\mathbb{E}\Big| \int_r^{\left( \lfloor \frac{r}{\delta _{\epsilon}} \rfloor +1 \right) \delta _{\epsilon}}{b_1\big( s\left( \delta _{\epsilon} \right) ,u_{s\left( \delta _{\epsilon} \right)}^{\epsilon},\hat{v}  _{r}^{\epsilon} \big) -\bar{b}_1\left( s\left( \delta _{\epsilon} \right) ,u_{s\left( \delta _{\epsilon} \right)}^{\epsilon} \right) ds} \Big|^2dr}\cr
	&\quad+C_{\alpha ,T}\int_0^{\left( \lfloor \frac{t}{\delta _{\epsilon}} \rfloor -3 \right) \delta _{\epsilon}}{\left( t-r \right) ^{-\alpha -\frac{3}{2}}\mathbb{E}\Big| \int_{\left( \lfloor \frac{r}{\delta _{\epsilon}} \rfloor +1 \right) \delta _{\epsilon}}^{\lfloor \frac{t}{\delta _{\epsilon}} \rfloor \delta _{\epsilon}}{b_1\big( s\left( \delta _{\epsilon} \right) ,u_{s\left( \delta _{\epsilon} \right)}^{\epsilon},\hat{v}  _{s}^{\epsilon} \big) -\bar{b}_1\left( s\left( \delta _{\epsilon} \right) ,u_{s\left( \delta _{\epsilon} \right)}^{\epsilon} \right) ds} \Big|^2dr}\cr
	&\quad+C_{\alpha ,T}\int_0^{\left( \lfloor \frac{t}{\delta _{\epsilon}} \rfloor -3 \right) \delta _{\epsilon}}{\left( t-r \right) ^{-\alpha -\frac{3}{2}}\mathbb{E}\Big| \int_{\lfloor \frac{t}{\delta _{\epsilon}} \rfloor \delta _{\epsilon}}^t{b_1\big( s\left( \delta _{\epsilon} \right) ,u_{s\left( \delta _{\epsilon} \right)}^{\epsilon},\hat{v}  _{s}^{\epsilon} \big) -\overline{b}_1\left( s\left( \delta _{\epsilon} \right) ,u_{s\left( \delta _{\epsilon} \right)}^{\epsilon} \right) ds} \Big|^2dr}\cr
	&\le C_{\alpha ,T}\int_0^{\left( \lfloor \frac{t}{\delta _{\epsilon}} \rfloor -3 \right) \delta _{\epsilon}}{\big( t-\big( \lfloor  {t}/{\delta _{\epsilon}} \rfloor -3 \big) \delta _{\epsilon} \big) ^{-\alpha -\frac{3}{2}}\big( \big( \lfloor  {r}/{\delta _{\epsilon}} \rfloor +1 \big) \delta _{\epsilon}-r \big)}\cr
	&\qquad\qquad\qquad\qquad\times \int_r^{\left( \lfloor \frac{r}{\delta _{\epsilon}} \rfloor +1 \right) \delta _{\epsilon}}{\mathbb{E}\big| b_1\big( s\left( \delta _{\epsilon} \right) ,u_{s\left( \delta _{\epsilon} \right)}^{\epsilon},\hat{v}  _{s}^{\epsilon} \big) -\bar{b}_1\left( s\left( \delta _{\epsilon} \right) ,u_{s\left( \delta _{\epsilon} \right)}^{\epsilon} \right) \big|^2ds}dr\cr
	&\quad+C_{\alpha ,T}\int_0^{\left( \lfloor \frac{t}{\delta _{\epsilon}} \rfloor -3 \right) \delta _{\epsilon}}{\left( t-r \right) ^{-\alpha -\frac{3}{2}}\left( \lfloor  {t}/{\delta _{\epsilon}} \rfloor -\lfloor {r}/{\delta _{\epsilon}} \rfloor -1 \right)}\cr
	&\qquad\qquad\qquad\qquad\times \sum_{k=\lfloor \frac{r}{\delta _{\epsilon}} \rfloor +1}^{\lfloor \frac{t}{\delta _{\epsilon}} \rfloor -1}{\mathbb{E}\Big| \int_{k\delta _{\epsilon}}^{\left( k+1 \right) \delta _{\epsilon}}{b_1\big( k\delta _{\epsilon},u_{k\delta _{\epsilon}}^{\epsilon},\hat{v}  _{s}^{\epsilon} \big) -\bar{b}_1\left( k\delta _{\epsilon},u_{k\delta _{\epsilon}}^{\epsilon} \right) ds} \Big|^2}dr\cr
	&\quad+C_{\alpha ,T}\int_0^{\left( \lfloor \frac{t}{\delta _{\epsilon}} \rfloor -3 \right) \delta _{\epsilon}}{\left( t-\left( \lfloor  {t}/{\delta _{\epsilon}} \rfloor -3 \right) \delta _{\epsilon} \right) ^{-\alpha -\frac{3}{2}}\left( t-\lfloor  {t}/{\delta _{\epsilon}} \rfloor \delta _{\epsilon} \right)}\cr &\qquad\qquad\qquad\qquad\times\int_{\lfloor \frac{t}{\delta _{\epsilon}} \rfloor \delta _{\epsilon}}^t{\mathbb{E}\left| b_1\big( s\left( \delta _{\epsilon} \right) ,u_{s\left( \delta _{\epsilon} \right)}^{\epsilon},\hat{v}  _{s}^{\epsilon} \big) -\bar{b}_1\left( s\left( \delta _{\epsilon} \right) ,u_{s\left( \delta _{\epsilon} \right)}^{\epsilon} \right) \right|^2ds}dr\cr
	&\le C_{\alpha ,x,y, T}\int_0^{\left( \lfloor \frac{t}{\delta _{\epsilon}} \rfloor -3 \right) \delta _{\epsilon}}{\left( 4\delta _{\epsilon} \right) ^{-\alpha -\frac{3}{2}}\delta _{\epsilon}^{2}dr}+C_{\alpha ,T}\int_0^{\left( \lfloor \frac{t}{\delta _{\epsilon}} \rfloor -3 \right) \delta _{\epsilon}}{\left( t-r \right) ^{-\alpha -\frac{3}{2}}}\left( \lfloor  {t}/{\delta _{\epsilon}} \rfloor -\lfloor  {r}/{\delta _{\epsilon}} \rfloor -1 \right) ^2\cr
	&\qquad\qquad\qquad\times\underset{\lfloor \frac{r}{\delta _{\epsilon}} \rfloor +1\le k\le \lfloor \frac{t}{\delta _{\epsilon}} \rfloor -1}{\max}\mathbb{E}\Big| \int_{k\delta _{\epsilon}}^{\left( k+1 \right) \delta _{\epsilon}}{b_1\big( k\delta _{\epsilon},u_{k\delta _{\epsilon}}^{\epsilon},\hat{v}  _{s}^{\epsilon} \big) -\bar{b}_1\left( k\delta _{\epsilon},u_{k\delta _{\epsilon}}^{\epsilon} \right) ds} \Big|^2dr\cr
	&\le \frac{C_{\alpha ,T}}{\delta _{\epsilon}^{2}}\int_0^t{\left( t-r \right) ^{\frac{1}{2}-\alpha}dr}\underset{0\le k\le \lfloor \frac{t}{\delta _{\epsilon}} \rfloor -1}{\max}\mathbb{E}\Big| \int_{k\delta _{\epsilon}}^{\left( k+1 \right) \delta _{\epsilon}}{b_1\big( k\delta _{\epsilon},u_{k\delta _{\epsilon}}^{\epsilon},\hat{v}  _{s}^{\epsilon} \big) -\bar{b}_1\left( k\delta _{\epsilon},u_{k\delta _{\epsilon}}^{\epsilon} \right) ds} \Big|^2\cr
	&\quad+C_{\alpha ,x,y, T}\delta _{\epsilon}^{\frac{1}{2}-\alpha}\cr
	&\le C_{\alpha ,x,y, T}\big( \delta _{\epsilon}^{\frac{1}{2}-\alpha}+ {\epsilon}/{\delta _{\epsilon}}\big) .\nonumber
	\end{align}
	For $ \mathcal{K}_{t}^{32}, $ using the assumption (A2)  and thanks to (\ref{en22}) and (\ref{en24}), it yields  
	\begin{align}
	\mathcal{K}_{t}^{32}
	%&=C_{\alpha ,T}\int_{\left( \lfloor \frac{t}{\delta _{\epsilon}} \rfloor -3 \right) \delta _{\epsilon}}^{\lfloor \frac{t}{\delta _{\epsilon}} \rfloor \delta _{\epsilon}}{\left( t-r \right) ^{-\alpha -\frac{3}{2}}\mathbb{E}\Big| \int_r^t{b_1\big( s\left( \delta _{\epsilon} \right) ,u_{s\left( \delta _{\epsilon} \right)}^{\epsilon},\hat{v}  _{s}^{\epsilon} \big) -\bar{b}_1\left( s\left( \delta _{\epsilon} \right) ,u_{s\left( \delta _{\epsilon} \right)}^{\epsilon} \right) ds} \Big|^2dr}\cr
	&\le C_{\alpha ,T}\int_{\left( \lfloor \frac{t}{\delta _{\epsilon}} \rfloor -3 \right) \delta _{\epsilon}}^{t}{\left( t-r \right) ^{-\alpha -\frac{1}{2}}\int_r^t{\mathbb{E}\big| b_1\big( s\left( \delta _{\epsilon} \right) ,u_{s\left( \delta _{\epsilon} \right)}^{\epsilon},\hat{v}  _{s}^{\epsilon} \big) -\bar{b}_1\left( s\left( \delta _{\epsilon} \right) ,u_{s\left( \delta _{\epsilon} \right)}^{\epsilon} \right) \big|^2ds}dr}\cr
	&\le C_{\alpha ,x,y, T}\int_{\left( \lfloor \frac{t}{\delta _{\epsilon}} \rfloor -3 \right) \delta _{\epsilon}}^{t}{\left( t-r \right) ^{\frac{1}{2}-\alpha}dr} 
	%&\le C_{\alpha ,x,y,T}\left( t-r \right) ^{\frac{3}{2}-\alpha}\left| _{\lfloor \frac{t}{\delta _{\epsilon}} \rfloor \delta _{\epsilon}}^{\left( \lfloor \frac{t}{\delta _{\epsilon}} \rfloor -3 \right) \delta _{\epsilon}} \right. \cr
	\le C_{\alpha ,x,y, T}\delta _{\epsilon}^{\frac{3}{2}-\alpha}.\nonumber
	\end{align}
	%Similarly, it is easy to prove that the following conclusions hold for $ \mathcal{K}_{t}^{33}  $
	%\begin{align} 
	%\mathcal{K}_{t}^{33}\le C_{\alpha ,x,y,x_0,t_1,t_2,T}\delta _{\epsilon}^{\frac{3}{2}-\alpha}.\nonumber
	%\end{align}
	Hence
	\begin{align}\label{en27}
	\mathcal{K}_{t}^{3}\leq C_{\alpha ,x,y, T}\big( \delta _{\epsilon}^{\frac{1}{2}-\alpha}+ {\epsilon}/{\delta _{\epsilon}}+\delta _{\epsilon}^{\frac{3}{2}-\alpha}\big).
	\end{align}
	Substituting (\ref{en26})-(\ref{en27}) into (\ref{en25}), it follows that (\ref{en28}) holds.
	%\begin{align}
	%\mathbb{E}\Big\|  \int_0^t{b_1\big( r\left( \delta _{\epsilon} \right) ,u_{r\left( \delta _{\epsilon} \right)}^{\epsilon},\hat{v}  _{r}^{\epsilon} \big) -\bar{b}_1\left( r\left( \delta _{\epsilon} \right) ,u_{r\left( \delta _{\epsilon} \right)}^{\epsilon} \right) dr} \Big\| _{\alpha}^{2}&\leq C_{\alpha ,x,y, T}\big(\frac{\epsilon}{\delta _{\epsilon}}+\delta _{\epsilon}^{2}+\delta _{\epsilon}^{\frac{1}{2}-\alpha}+ \delta _{\epsilon}^{\frac{3}{2}-\alpha}\big).\nonumber
	%\end{align}
	The proof is completed.\qed
	
	Thanks to (\ref{en29}), (\ref{en30}) and \lemref{lem6.2}, it follows that
	\begin{align}\label{en31}
	\mathcal{I}_{t}^{1}&\leq  C_{\alpha ,R,T}\int_0^t{\left( t-r \right) ^{-2\alpha}\mathbb{E}\left[ \left\|  u_{r}^{\epsilon}-\bar{u}_r \right\|  _{\alpha}^{2}\mathbf{1}_{\varOmega _M\cap \{T\leq \tau _{R}^{\epsilon}\}} \right] dr}\cr
	&\quad+C_{\alpha ,x,y,R,T}\big( \delta _{\epsilon}^{2\kappa \land 2\iota \land \left( {1}/{2}-\alpha\right) }+ {\epsilon}/{\delta _{\epsilon}}\big) .
	\end{align}
	Then, substituting (\ref{en32})  and (\ref{en31}) into (\ref{en34}), it yields
	\begin{align}
	\mathbb{E}\left[ \left\|  u_{t}^{\epsilon}-\bar{u}_t \right\|  _{\alpha}^{2}\mathbf{1}_{\varOmega _M\cap \{T\leq \tau _{R}^{\epsilon}\}} \right]
	&\leq C_{\alpha ,R,M,T} \int_0^t{   \left( t-r \right) ^{-\frac{1}{2}-\alpha} r^{-\frac{1}{2}-\alpha}   \mathbb{E}\left[ \left\|  u_{r}^{\epsilon}-\bar{u}_r \right\|  _{\alpha}^{2}\mathbf{1}_{\varOmega _M\cap \{T\leq \tau _{R}^{\epsilon}\}} \right] dr}\cr
	& \quad+C_{\alpha ,R,x,y, T}\big( \delta _{\epsilon}^{2\kappa \land 2\iota \land \left( {1}/{2}-\alpha\right) }+ {\epsilon}/{\delta _{\epsilon}}\big).\nonumber
	\end{align}
	According to the Gronwall-type lemma (Lemma 7.6 in  \cite{Nualart2002Differential}), we deduce (\ref{en49}).
	%\begin{align}
	%\mathbb{E}\left[ \left\|  u_{t}^{\epsilon}-\bar{u}_t \right\|  _{\alpha}^{2}\mathbf{1}_{\varOmega _M\cap \{T\leq \tau _{R}^{\epsilon}\}} \right]
	%%&\leq C_{\alpha ,x,y,T}\big(\delta _{\epsilon}^{2\kappa \land 2\iota \land \left( {1}/{2}-\alpha\right) }+ {\epsilon}/{\delta _{\epsilon}}\big)d_\alpha\text{exp} \big[c_\alpha t \big(C_{\alpha ,M,T}\left( 1+\varLambda _{\alpha}^{2}\left( B^H \right)\right)^{\frac{1}{ 1/2-\alpha} }\big)\big]\cr
	%&\leq C_{\alpha ,x,y,x_0,t_1,t_2,R,M,T}\big(\delta _{\epsilon}^{2\kappa \land 2\iota \land \left( {1}/{2} -\alpha\right)}+ {\epsilon}/{\delta _{\epsilon}}\big).\nonumber
	%\end{align}
	The proof is completed.\qed

	Finally, the proof of our main result can be finished. 
	\para{Proof of \thmref{thm2.1}:}
	Selecting  $ \delta _{\epsilon}=\epsilon \ln ^{\epsilon ^{-\varsigma}} \left(    0<\varsigma<1\right) , $ we can get
	\begin{align}
	\underset{\epsilon \rightarrow 0}{\lim}\underset{t\in \left[ 0,T \right]}{\text{sup}}\mathbb{E}\left[ \left\|  u_{t}^{\epsilon}-\bar{u}_t \right\|  _{\alpha}^{2}\mathbf{1}_{\varOmega _M\cap \{T\leq \tau _{R}^{\epsilon}\}} \right]=0.\nonumber
	\end{align}
	Then, let $ M\rightarrow+\infty, $ as $ \Omega_M\nearrow\Omega, $ it yields
	\begin{align}\label{en47}
	\underset{\epsilon \rightarrow 0}{\lim}\underset{t\in \left[ 0,T \right]}{\text{sup}}\mathbb{E}\left[ \left\|  u_{t}^{\epsilon}-\bar{u}_t \right\|  _{\alpha}^{2}\mathbf{1}_{ \{T\leq \tau _{R}^{\epsilon}\}} \right]=0. 
	\end{align}
	
	Using the H\"{o}lder inequality and Chebushev's inequality, thanks to (\ref{en46}), (\ref{en22}) and Lemma 7.5 in \cite{Nualart2002Differential}, it yields
	\begin{align}\label{en48}
	&\quad\ \underset{t\in \left[ 0,T \right]}{\text{sup}}\mathbb{E}\left[ \left\|  u_{t}^{\epsilon}-\bar{u}_t \right\|  _{\alpha}^{2}\mathbf{1}_{\{T>\tau _{R}^{\epsilon}\}} \right] \le \underset{t\in \left[ 0,T \right]}{\text{sup}}\left[ \mathbb{E}\left\|  u_{t}^{\epsilon}-\bar{u}_t \right\|  _{\alpha}^{4} \right] ^{\frac{1}{2}}\left[ \mathbb{P}\left( T>\tau _{R}^{\epsilon} \right) \right] ^{\frac{1}{2}}\cr
	&\le \underset{t\in \left[ 0,T \right]}{\text{sup}}R^{-1}\left[ \mathbb{E}\left\|  u_{t}^{\epsilon} \right\|  _{\alpha}^{4}+\mathbb{E}\left\|  \bar{u}_t \right\|  _{\alpha}^{4} \right] ^{\frac{1}{2}}\Big[ \mathbb{E}\Big(  \sup_{t\in \left[ 0,T \right]}\left\|  u_{t}^{\epsilon} \right\|  _{\alpha}^{2}+\sup_{t\in \left[ 0,T \right]}\left\|  \bar{u}_t \right\|  _{\alpha}^{2}+\sup_{t\in \left[ 0,T \right]}
	\left( \varLambda _{\alpha}^{0,t}\left( B^H \right) \right) ^2\Big)   \Big] ^{\frac{1}{2}}\cr
	&\le C_{\alpha ,x,  T}R^{-1}.
	\end{align}
	Note that
	\begin{align}
	\underset{t\in \left[ 0,T \right]}{\text{sup}}\mathbb{E}\left\|  u_{t}^{\epsilon}-\bar{u}_t \right\|  _{\alpha}^{2}=\underset{t\in \left[ 0,T \right]}{\text{sup}}\mathbb{E}\left[ \left\|  u_{t}^{\epsilon}-\bar{u}_t \right\|  _{\alpha}^{2}\mathbf{1}_{\{T\leq \tau _{R}^{\epsilon}\}} \right] +\underset{t\in \left[ 0,T \right]}{\text{sup}}\mathbb{E}\left[ \left\|  u_{t}^{\epsilon}-\bar{u}_t \right\|  _{\alpha}^{2}\mathbf{1}_{\{T>\tau _{R}^{\epsilon}\}} \right].
	\end{align}
	Thanks to (\ref{en47}) and (\ref{en48}), let $ \epsilon \rightarrow 0 $ firstly and $ R\rightarrow \infty $ secondly, we can get the desired estimate (\ref{en35}). This completes the proof of \thmref{thm2.1}.\qed

	\section*{Acknowledgments}
	This work was partly supported by the National Natural Science Foundation of China under Grant Nos. 11772255 and 12072264, the National Key Research and Development Program of China under Grant No. 2018AAA0102201, the Fundamental Research Funds for the Central Universities, the Research Funds for Interdisciplinary Subject of Northwestern Polytechnical University, the Shaanxi Project for Distinguished Young Scholars, the Shaanxi Provincial Key R\&D Program 2020KW-013 and 2019TD-010. 
	
	\section*{Appendix}
	In this section, we give the detailed proofs of \lemref{lem4.1} and  \lemref{lemm6.1}:
	\para{Proof of \lemref{lem4.1}:} According to the estimate (\ref{en13}) and the classical Bogoliubov-Krylov argument, it is possible to get that the existence of an invariant measure $ \mu ^{s,x}.  $  Then, for any Lipschitz function $ \varphi:\mathbb{R}^m\rightarrow\mathbb{R}, $ 
	thanks to (\ref{en14}) and (\ref{en15}), we obtain
	\begin{align}
	\Big| P_{t}^{s,x}\varphi \left( y \right) -\int_{\mathbb{R}^m}{\varphi \left( z \right) \mu ^{s,x}\left( dz \right)} \Big|
	%&=\Big| P_{t}^{s,x}\varphi \left( y \right) -\int_{\mathbb{R}^m}{P_{t}^{s,x}\varphi \left( z \right) \mu ^{s,x}\left( dz \right)} \Big|\cr
	%&=\Big| \int_{\mathbb{R}^m}{\left( \mathbb{E}\varphi \left( Y_{t}^{s,x,y} \right) -\mathbb{E}\varphi \left( Y_{t}^{s,x,z} \right) \right) \mu ^{s,x}\left( dz \right)} \Big|\cr
	&\le \int_{\mathbb{R}^m}{\left| P_{t}^{s,x}\varphi \left( y \right) -P_{t}^{s,x}\varphi \left( z \right) \right|\mu ^{s,x}\left( dz \right)}\cr
	%&\le \int_{\mathbb{R}^m}{\left| \mathbb{E}\varphi \left( Y_{t}^{s,x,y} \right) -\mathbb{E}\varphi \left( Y_{t}^{s,x,z} \right) \right|\mu ^{s,x}\left( dz \right)}\cr
	&\le \int_{\mathbb{R}^m}{\mathbb{E}\left| \varphi \left( v_{t}^{s,x,y} \right) -\varphi \left( v_{t}^{s,x,z} \right) \right|\mu ^{s,x}\left( dz \right)}\cr
	&\le L_{\varphi}\int_{\mathbb{R}^m}{\mathbb{E}\left| v_{t}^{s,x,y}-v_{t}^{s,x,z} \right|\mu ^{s,x}\left( dz \right)}\cr
	%&\le L_{\varphi}\int_{\mathbb{R}^m}{\left[ \mathbb{E}\left| Y_{t}^{s,x,y}-Y_{t}^{s,x,z} \right|^2 \right] ^{\frac{1}{2}}\mu ^{s,x}\left( dz \right)}\cr
	&\le L_{\varphi}\int_{\mathbb{R}^m}{e^{-\frac{\beta _1}{2}t}\left| y-z \right|\mu ^{s,x}\left( dz \right)}\cr
	&\le C L_{\varphi}e^{-\frac{\beta _1}{2}t}\left( 1+\left| x \right|+\left| y \right| \right),\nonumber
	\end{align}
	where $ L_{\varphi}=\underset{x\ne y}{\text{sup}}\frac{\left| \varphi \left( x \right) -\varphi \left( y \right) \right|}{\left| x-y \right|}. $ The  invariant measure $ \mu ^{s,x} $ is unique and strong mixing.
	
	Thanks to (\ref{en13}), for all $ t>0, $ we have
	\begin{align}
	\int_{\mathbb{R}^m}{\left| z \right|^p\mu ^{s,x}\left( dz \right)}& =\int_{\mathbb{R}^m}{\mathbb{E}\left| v_{t}^{s,x,z} \right|^p\mu ^{s,x}\left( dz \right)}\le C_{p,T}\left( 1+\left| x \right|^p \right) +\int_{\mathbb{R}^m}{e^{-\beta _{p}^{*} t}\left| z \right|^p}\mu ^{s,x}\left( dz \right).\nonumber
	\end{align}
	Therefore, if we take $ t>0 $ such that $ e^{-\beta _{p}^{*}t} \leq 1/2, $ we can get (\ref{en15}).
	Moreover, due to (\ref{en13}), (\ref{en14}) and (\ref{en15}), we also have
	\begin{align}
	&\qquad\Big| \mathbb{E}b_1\left( s,x,v_{t}^{s,x,y} \right) -\int_{\mathbb{R}^m}{b_1\left( s,x,z \right) \mu ^{s,x}\left( dz \right)} \Big|\cr
	&=\Big| \int_{\mathbb{R}^m}{\mathbb{E}b_1\left( s,x,v_{t}^{s,x,y} \right) -\mathbb{E}b_1\left( s,x,v_{t}^{s,x,z} \right) \mu ^{s,x}\left( dz \right)} \Big|\cr
	&\le C_T\int_{\mathbb{R}^m}{\mathbb{E}\big[ \left| v_{t}^{s,x,y}-v_{t}^{s,x,z} \right|\big( 1+\left| x \right|^{\theta _2}+\left| v_{t}^{s,x,y} \right|^{\theta _3}+\left| v_{t}^{s,x,z} \right|^{\theta _3} \big) \big]}\mu ^{s,x}\left( dz \right) \cr
	&\le C_T\int_{\mathbb{R}^m}{\big[ \mathbb{E}\left| v_{t}^{s,x,y}-v_{t}^{s,x,z} \right|^2\mathbb{E}\big( 1+\left| x \right|^{2\theta _2}+\left| v_{t}^{s,x,y} \right|^{2\theta _3}+\left| v_{t}^{s,x,z} \right|^{2\theta _3} \big) \big] ^{\frac{1}{2}}}\mu ^{s,x}\left( dz \right) \cr
	&\le C_Te^{-\frac{\beta _1}{2}t}\int_{\mathbb{R}^m}{\left| z-y \right|\big( 1+\left| x \right|^{\theta _2}+\left| x \right|^{\theta _3}+\left| y \right|^{\theta _3}+\left| z \right|^{\theta _3} \big)}\mu ^{s,x}\left( dz \right) \cr
	&\le C_Te^{-\frac{\beta _1}{2}t}\int_{\mathbb{R}^m}{\big( 1+\left| x \right|^{\theta _2\lor \theta _3}+\left| y \right|^{\theta _3\lor 1}+\left| z \right|^{\theta _3\lor 1} \big) ^2}\mu ^{s,x}\left( dz \right) \cr
	&\le C_Te^{-\frac{\beta _1}{2}t}\big( 1+\left| x \right|^{2\left( \theta _2\lor \theta _3\lor 1 \right)}+\left| y \right|^{2\left( \theta _3\lor 1 \right)} \big) .\nonumber
	\end{align}
	The proof is completed.\qed
	\para{Proof of \lemref{lemm6.1}:} For any $ x_1,x_2\in \mathbb{R}^n   $ with  $ \left|x_i \right|\leq R, $ thanks to the assumptions and equation (\ref{en13}), (\ref{en17})
	and (\ref{en16}),  we can get
	\begin{align}
	\left| \bar{b}_1\left( s_1,x_1 \right) -\bar{b}_1\left( s_2,x_2 \right) \right|&\le \Big| \int_{\mathbb{R}^m}{b_1\left( s_1,x_1,z \right) \mu ^{s_1,x_1}\left( dz \right)}-\mathbb{E}b_1\left( s_1,x_1,v_{t}^{s_1,x_1,0} \right) \Big|\cr
	&\quad+\left| \mathbb{E}b_1\left( s_1,x_1,v_{t}^{s_1,x_1,0} \right) -\mathbb{E}b_1\left( s_2,x_1,v_{t}^{s_1,x_1,0} \right) \right|\cr
	&\quad+\left| \mathbb{E}b_1\left( s_2,x_1,v_{t}^{s_1,x_1,0} \right) -\mathbb{E}b_1\left( s_2,x_2,v_{t}^{s_1,x_1,0} \right) \right|\cr
	&\quad+\left| \mathbb{E}b_1\left( s_2,x_2,v_{t}^{s_1,x_1,0} \right) -\mathbb{E}b_1\left( s_2,x_2,v_{t}^{s_2,x_2,0} \right) \right|\cr
	&\quad+\Big| \mathbb{E}b_1\left( s_2,x_2,v_{t}^{s_2,x_2,0} \right) -\int_{\mathbb{R}^m}{b_1\left( s_2,x_2,z \right) \mu ^{s_2,x_2}\left( dz \right)} \Big|\cr
	&\le C_Te^{-\frac{\beta _1}{2}t}\big( 1+\left| x_1 \right|^{2\left( \theta _2\lor \theta _3\lor 1 \right)}+\left| x_2 \right|^{2\left( \theta _2\lor \theta _3\lor 1 \right)} \big) \cr
	&\quad+C_T\left| s_1-s_2 \right|^{\kappa}\big( 1+\left| x_1 \right|^{\theta _2}+\mathbb{E}\left| v_{t}^{s_1,x_1,0} \right|^{\theta _3} \big) \cr
	&\quad+C_{R,T}\left| x_1-x_2 \right|\big( 1+\mathbb{E}\left| v_{t}^{s_1,x_1,0} \right|^{\theta _1} \big) \cr
	&\quad+C_T\mathbb{E}\big[ \left| v_{t}^{s_1,x_1,0}-v_{t}^{s_2,x_2,0} \right|\big( 1+\left| x_2 \right|^{\theta _2}+\left| v_{t}^{s_1,x_1,0} \right|^{\theta _3}+\left| v_{t}^{s_2,x_2,0} \right|^{\theta _3} \big) \big] \cr
	&\le C_Te^{-\frac{\beta _1}{2}t}\big( 1+\left| x_1 \right|^{2\left( \theta _2\lor \theta _3\lor 1 \right)}+\left| x_2 \right|^{2\left( \theta _2\lor \theta _3\lor 1 \right)} \big) \cr
	&\quad+C_T\left| s_1-s_2 \right|^{\kappa}\big( 1+\left| x_1 \right|^{\theta _2}+\left| x_1 \right|^{\theta _3} \big) +C_{R,T}\left| x_1-x_2 \right|\big( 1+\left| x_1 \right|^{\theta _1} \big) \cr
	&\quad+C_T\left( \left| s_1-s_2 \right|^{\iota}+\left| x_1-x_2 \right| \right) \big( 1+\left| x_1 \right|^{2\left( \alpha _1\lor \theta _3 \right)}+\left| x_2 \right|^{2\left( \alpha _1\lor \alpha _2\lor \theta _2\lor \theta _3 \right)} \big) \cr
	&\le C_{R,T}e^{-\frac{\beta _1}{2}t} +C_{R,T}\left( \left| s_1-s_2 \right|^{\kappa}+\left| s_1-s_2 \right|^{\iota}+\left| x_1-x_2 \right| \right).\nonumber
	%&\qquad\qquad\quad\times \big( 1+\left| x_1 \right|^{\theta _1\lor \theta _2\lor 2\left( \alpha _1\lor \theta _3 \right)}+\left| x_2 \right|^{2\left( \alpha _1\lor \alpha _2\lor \theta _2\lor \theta _3 \right)} \big) .\nonumber
	\end{align}
	Let $ t\rightarrow +\infty, $ we obtain (\ref{en19}). Moreover, thanks to (\ref{en15}), we also have
	\begin{align}
	\left| \bar{b}_1\left( t,x \right) \right|
	%&=\Big| \int_{\mathbb{R}^m}{b_1\left( t,x,z \right) \mu ^{t,x}\left( dz \right)} \Big|
	\le \int_{\mathbb{R}^m}{\left| b_1\left( t,x,z \right) \right|\mu ^{t,x}\left( dz \right)} 
	\le C_T\int_{\mathbb{R}^m}{\left( 1+\left| x \right|+\left| z \right| \right) \mu ^{t,x}\left( dz \right)}\le C_T\left( 1+\left| x \right| \right).\nonumber
	\end{align}
	The proof is completed.\qed

	\bibliography{references}
\end{document}